\documentclass[final,3p]{elsarticle}
 \usepackage{graphics}
 \usepackage{graphicx}
 \usepackage{epsfig}
\usepackage{amssymb}
 \usepackage{amsthm}
 \usepackage{lineno}
 \usepackage{amsmath}
   \numberwithin{equation}{section}
\usepackage{mathrsfs}

\newtheorem{thm}{Theorem}[section]

\newtheorem{lem}[thm]{Lemma}

\biboptions{sort&compress,square}
\allowdisplaybreaks
\begin{document}
\begin{frontmatter}
\author{Yuchen Yang}
\ead{yangyc580@nenu.edu.cn}
\author{Tong Wu\corref{cor3}}
\ead{wut977@nenu.edu.cn}
\cortext[cor3]{Corresponding author.}

\address{School of Mathematics and Statistics, Northeast Normal University,
Changchun, 130024, China}
\title{The Spectral Einstein functional and Kastler-Kalau-Walze type theorems}
\begin{abstract}
In this paper, on the basis of defining the spectral Einstein functional associated with the Dirac operator for manifolds with boundary, we prove Kastler-Kalau-Walze type theorem for the spectral Einstein functional associated with the Dirac operator on low-dimensional manifolds with boundary.
\end{abstract}
\begin{keyword}The spectral Einstein functional; the Dirac operator; Kastler-Kalau-Walze type theorem;\\

\end{keyword}
\end{frontmatter}
\textit{2010 Mathematics Subject Classification:}
53C40; 53C42.
\section{Introduction}
 In recent years, many geometers have made more in-depth study of noncommutative residues.
 In \cite{Gu}, Guillemin found noncommutative residues are of great importance to the study of noncommutative geometry.
 In the 1990s, the scholar Connes induced a conformal 4-dimensional Polyakov action analogy by using non-commutative residues\cite{Co1}, he also showed that the non-commutative residues on the compact manifold equals the Dixmier's trace on the pseudo-differential operators with order $-{\rm {dim}}M$ in \cite{Co2}, so that the non-commutative residues can be used as a non-commutative integral.
 Connes conjectured that the noncommutative residue of the $-2$ power of the Dirac operator is proportional to the Einstein-Hilbert action.
In 1995, Kastler\cite{Ka} proved this conjecture. At the same time, Kalau and Walze also proved this conjecture in the normal coordinate system in \cite{KW}.
In 1996, Ackermann  found that the noncommutative residue of the inverse square of the Dirac operator $\widetilde{{\rm  Wres}}(D^{-2})$ is essentially the second coefficient of the heat kernel expansion of the square of the Dirac operator in \cite{Ac}.
In the same year, Fedosov et al. proposed the idea of noncommutative residue on manifolds with boundary. They defined the non-commutative residue on Boutet de Monvel algebra and proved that it is a unique continuous trace in \cite{FGLS}.

On the other hand, Wang provides a kind of method to study the Kastler-Kalau-Walze type theorem for manifolds with boundary. In the early 21st century, Chinese scholar Wang Yong extended the results of Connes to manifolds with boundary in \cite{Wa1,Wa2}, studied the conformal invariant problem and gravity action of manifolds with boundary, and proved the  Kastler-Kalau-Walze type theorems of Dirac operator and the Signature operator on low dimensional manifolds with boundary \cite{Wa3}.
In \cite{Wa4}, Wang defined lower dimensional volumes of spin manifolds with boundary. And computed the lower dimensional volume $Vol^{(2,2)}_6 $ and got a Kastler-Kalau-Walze type theorem in 6-dimensional spin manifolds with boundary and the associated Dirac operator D.
They computed the lower dimensional volume $Vol^{(1,3)}_6$ and got a Kastler-Kalau-Walze type theorem in 6-dimensional spin manifolds with boundary and the associated Dirac operator $D$ and $D^3$ in \cite{Wa5}.
  In \cite{DL}, Dabrowski L, Sitarz A and Zalecki P defined bilinear functionals of vector fields and differential forms, the densities of which yield the metric and Einstein tensors on even-dimensional Riemannian manifolds and proved that for the conformally rescaled geometry of the noncommutative two-torus the Einstein functional vanishes.
 Motivated by \cite{DL}, we define  the spectral Einstein functional associated with the Dirac operator for manifolds with boundary, and the motivation of this paper is to compute the noncommutative residue
 $\widetilde{{\rm Wres}}[\pi^+(\nabla_X^{S(TM)}\nabla_Y^{S(TM)}D^{-1})\circ\pi^+(D^{-3})]$ on 4-dimensional spin manifolds with boundary and $\widetilde{{\rm Wres}}[\pi^+(\nabla_X^{S(TM)}\nabla_Y^{S(TM)}D^{-3})\circ\pi^+(D^{-3})]$ on 6-dimensional spin manifolds with boundary. \\
\indent The paper is organized in the following way. In Section \ref{section:2}, we define the spectral Einstein functional associated with the Dirac operator and and get the noncommutative residue for manifolds without boundary. In Section \ref{section:3}, we prove the Kastler-Kalau-Walze type theorem for the spectral Einstein functional associated with the Dirac operator on 4-dimensional manifolds with boundary and we get a Kastler-Kalau-Walze type theorem associated for 3-dimensional spin manifolds in Section \ref{section:4}. In Section \ref{section:5}, for 6-dimensional manifolds with boundary and the spectral Einstein functional associated with the Dirac operator, we compute the Kastler-Kalau-Walze type theorem.
\section{The spectral Einstein functional associated with the Dirac operator}
\label{section:2}
Firstly we review the relevant knowledge of the Dirac operator. We have $M$ is an $n$-dimensional ($n\geq 3$) oriented compact Riemannian manifold with a Riemannian metric $g^{M}$ and let $\nabla^L$ be the Levi-Civita connection about $g^{M}$. In the fixed orthonormal frame $\{e_1,\cdots,e_n\}$, the connection matrix $(\omega_{s,t})$ is defined by
\begin{equation}
\label{a2}
\nabla^L(e_1,\cdots,e_n)= (e_1,\cdots,e_n)(\omega_{s,t}).
\end{equation}
\indent Let $\epsilon (e_j^*)$,~$\iota (e_j^*)$ be the exterior and interior multiplications respectively, where $e_j^*=g^{TM}(e_j,\cdot)$.
Write
\begin{equation}
\label{a3}
\widehat{c}(e_j)=\epsilon (e_j^* )+\iota
(e_j^*);~~
c(e_j)=\epsilon (e_j^* )-\iota (e_j^* ),
\end{equation}
which satisfies
\begin{align}
\label{a4}
&\widehat{c}(e_i)\widehat{c}(e_j)+\widehat{c}(e_j)\widehat{c}(e_i)=2g^{M}(e_i,e_j);~~\nonumber\\
&c(e_i)c(e_j)+c(e_j)c(e_i)=-2g^{M}(e_i,e_j);~~\nonumber\\
&c(e_i)\widehat{c}(e_j)+\widehat{c}(e_j)c(e_i)=0.\nonumber\\
\end{align}
By \cite{Y}, we have the Dirac operator
\begin{align}
\label{a5}
&D=\sum^n_{i=1}c(e_i)[e_i-\frac{1}{4}\sum_{s,t}\omega_{s,t}
(e_i)c(e_s)c(e_t)].\nonumber\\
\end{align}
\indent We define spin connection by $\nabla_X^{S(TM)}:=X+\frac{1}{4}\Sigma_{ij}\langle\nabla_X^L{e_i},e_j\rangle c(e_i)c(e_j)$. Let $A(X)=\frac{1}{4}\Sigma_{ij}\langle\nabla_X^L{e_i},e_j\rangle c(e_i)c(e_j)$, so we have
\begin{align}\label{ddd}
\nabla_X^{S(TM)}\nabla_Y^{S(TM)}&=[X+A(X)][Y+A(Y)]\nonumber\\
&=XY+X\cdot A(Y)+A(X)Y+A(X)A(Y)\nonumber\\
&=XY+X[A(Y)]+A(Y)X+A(X)Y+A(X)A(Y),\nonumber\\
\end{align}
where $X=\Sigma_{j=1}^nX_j\partial_{x_j}, Y=\Sigma_{l=1}^nY_l\partial_{x_l}$ and  $XY=\Sigma_{j,l=1}^n(X_jY_l\partial_{x_j}\partial_{x_l}+X_j\frac{\partial_{Y_l}}{\partial_{x_j}}\partial_{x_l}).$ \\
\indent Let $g^{ij}=g(dx_{i},dx_{j})$, $\xi=\sum_{k}\xi_{j}dx_{j}$ and  $\nabla^L_{\partial_{i}}\partial_{j}=\sum_{k}\Gamma_{ij}^{k}\partial_{k}$,  we denote that
\begin{align}
&\sigma_{i}=-\frac{1}{4}\sum_{s,t}\omega_{s,t}
(e_i)c(e_s)c(e_t)
;~~~\xi^{j}=g^{ij}\xi_{i};~~~~\Gamma^{k}=g^{ij}\Gamma_{ij}^{k};~~~~\sigma^{j}=g^{ij}\sigma_{i};
~~~~a^{j}=g^{ij}a_{i}.
\end{align}
And by $\partial_{x_j}=-\sqrt{-1}\xi_j$, we have the following lemmas.
\begin{lem}\label{lem3} The following identities hold:
\begin{align}
\label{b22}
&\sigma_{0}(\nabla_X^{S(TM)}\nabla_Y^{S(TM)})=X[A(Y)]+A(X)A(Y);\nonumber\\
&\sigma_{1}(\nabla_X^{S(TM)}\nabla_Y^{S(TM)})=\sqrt{-1}\Sigma_{j,l=1}^nX_j\frac{\partial_{Y_l}}{\partial_{x_j}}\partial_{x_l}+\sqrt{-1}\Sigma_jA(Y)X_j\xi_j+\sqrt{-1}\Sigma_lA(Y)Y_l\xi_l;\nonumber\\
&\sigma_{2}(\nabla_X^{S(TM)}\nabla_Y^{S(TM)})=-\Sigma_{j,l=1}^nX_jY_l\xi_j\xi_l.\nonumber\\
\end{align}
\end{lem}
\begin{lem}
By\cite{Wa5}\label{lem356}, the following identities hold:
\begin{align}
\label{b22222}
&\sigma_{-1}(D^{-1})=\frac{\sqrt{-1}c(\xi)}{|\xi|^{2}};\nonumber\\
&\sigma_{-2}(D^{-1})=\frac{c(\xi)\sigma_0(D)c(\xi)}{|\xi|^4}
+\frac{c(\xi)}{|\xi|^6}\Sigma_jc(\mathrm{d}x_j)(\partial_{x_j}[c(\xi)]|\xi|^2-c(\xi)\partial_{x_j}(|\xi|^2));\nonumber\\
&\sigma_{-3}(D^{-3})=-\sqrt{-1}|\xi|^{-4}c(\xi);\nonumber\\
&\sigma_{-4}(D^{-3})=\frac{c(\xi)\sigma_2(D^3)c(\xi)}{|\xi|^8}+\frac{\sqrt{-1}c(\xi)}{|\xi|^8}\bigg(|\xi|^4c(dx_n)\partial_{x_n}c(\xi')
                     -2h'(0)c(\mathrm{d}x_n)c(\xi)\nonumber\\
                   &~~~~~+2\xi_nc(\xi)\partial{x_n}c(\xi')+4\xi_nh'(0)\bigg);\nonumber\\
\end{align}
\end{lem}
\indent We write
 \begin{eqnarray}
D_x^{\alpha}&=(-i)^{|\alpha|}\partial_x^{\alpha};
~\sigma(D_t^3)=p_3+p_2+p_1+p_0;
~(\sigma(D_t)^{-3})=\sum^{\infty}_{j=3}q_{-j}.
\end{eqnarray}

\indent By the composition formula of pseudodifferential operators, we have
\begin{align}
1=\sigma(D^3\circ D^{-3})
&=\sum_{\alpha}\frac{1}{\alpha!}\partial^{\alpha}_{\xi}[\sigma(D]
D_x^{\alpha}[\sigma(D^{-3})]\nonumber\\
&=(p_3+p_2+p_1+p_0)(q_{-3}+q_{-4}+q_{-5}+\cdots)\nonumber\\
&~~~+\sum_j(\partial_{\xi_j}p_3+\partial_{\xi_j}p_2++\partial_{\xi_j}p_1+\partial_{\xi_j}p_0)
(D_{x_j}q_{-3}+D_{x_j}q_{-4}+D_{x_j}q_{-5}+\cdots)\nonumber\\
&=p_3q_{-3}+(p_3q_{-4}+p_2q_{-3}+\sum_j\partial_{\xi_j}p_3D_{x_j}q_{-3})+\cdots,
\end{align}
so
\begin{equation}
q_{-3}=p_3^{-1};~q_{-4}=-p_3^{-1}[p_2p_3^{-1}+\sum_j\partial_{\xi_j}p_3D_{x_j}(p_-3^{-1})].
\end{equation}
Then,
\begin{align}\label{mki}
\sigma_{1}(\nabla_X^{S(TM)}\nabla_Y^{S(TM)}D^{-1})&=-\sqrt{-1}\Sigma_{j,l=1}^nX_jY_l\xi_j\xi_lc(\xi)|\xi|^{-2};\nonumber\\
\sigma_{-1}(\nabla_X^{S(TM)}\nabla_Y^{S(TM)}D^{-3})&=\sqrt{-1}\Sigma_{j,l=1}^nX_jY_l\xi_j\xi_lc(\xi)|\xi|^{-4}.
\end{align}

According to \cite{DL}, we have the following theorem

\begin{thm}\label{thm2} If $M$ is an $n$-dimensional compact oriented manifolds without boundary, and $n$ is even, then we get the following equality:
\begin{align}
\label{a29}
{\rm Wres}[\nabla_X^{S(TM)}\nabla_Y^{S(TM)}D^{-n}]
&=\frac{(2\pi)^{\frac{n}{2}}}{3(\frac{n}{2}-1)!}\int_{M}[Ric(X,Y)-\frac{1}{2}sg(X,Y)]d{\rm Vol_{M}}+\frac{(2\pi)^{\frac{n}{2}}}{4(\frac{n}{2}-1)!}\int_{M}sg(X,Y)d{\rm Vol_{M}},
\end{align}
where $s$ is the scalar curvature, ${\rm Vol_{M}}$ is the volume of $M$ and $Ric$ denotes Ricci tensor on $M$.
\end{thm}

\section{The noncommutative residue for $4$-dimensional manifolds with boundary}
\label{section:3}
 In this section, we compute the noncommutative residue $\widetilde{{\rm Wres}}[\pi^+(\nabla_X^{S(TM)}\nabla_Y^{S(TM)}D^{-1})\circ\pi^+(D^{-3})]$ on $4$-dimensional oriented compact manifolds with boundary. We review some basic facts and formulas about Boutet de Monvel's calculus and the definition of the noncommutative residue for manifolds with boundary which will be used in the following. Details are detailed in the the Section 2 in \cite{Wa3}.\\
 \indent Let $M$ be a 4-dimensional compact oriented manifold with boundary $\partial M$.
We assume that the metric $g^{M}$ on $M$ has the following form near the boundary,
\begin{equation}
\label{b1}
g^{M}=\frac{1}{h(x_{n})}g^{\partial M}+dx _{n}^{2},
\end{equation}
where $g^{\partial M}$ is the metric on $\partial M$ and $h(x_n)\in C^{\infty}([0, 1)):=\{\widehat{h}|_{[0,1)}|\widehat{h}\in C^{\infty}((-\varepsilon,1))\}$ for
some $\varepsilon>0$ and $h(x_n)$ satisfies $h(x_n)>0$, $h(0)=1,$ where $x_n$ denotes the normal directional coordinate.\\
\indent Similar to \cite{Wa3},  then we can compute the noncommutative residue
\begin{align}
\label{b14}
&\widetilde{{\rm Wres}}[\pi^+(\nabla_X^{S(TM)}\nabla_Y^{S(TM)}D^{-1})\circ\pi^+(D^{-3})]\nonumber\\
&=\int_M\int_{|\xi|=1}{\rm
trace}_{\wedge^*T^*M\bigotimes\mathbb{C}}[\sigma_{-4}(\nabla_X^{S(TM)}\nabla_Y^{S(TM)}D^{-1}\circ D^{-3})]\sigma(\xi)dx+\int_{\partial M}\Phi,
\end{align}
where
\begin{align}
\label{b15}
\Phi &=\int_{|\xi'|=1}\int^{+\infty}_{-\infty}\sum^{\infty}_{j, k=0}\sum\frac{(-i)^{|\alpha|+j+k+1}}{\alpha!(j+k+1)!}
\times {\rm trace}_{\wedge^*T^*M\bigotimes\mathbb{C}}[\partial^j_{x_n}\partial^\alpha_{\xi'}\partial^k_{\xi_n}\sigma^+_{r}(\nabla_X^{S(TM)}\nabla_Y^{S(TM)}D^{-1})(x',0,\xi',\xi_n)
\nonumber\\
&\times\partial^\alpha_{x'}\partial^{j+1}_{\xi_n}\partial^k_{x_n}\sigma_{l}(D^{-3})(x',0,\xi',\xi_n)]d\xi_n\sigma(\xi')dx',
\end{align}
and the sum is taken over $r+l-k-j-|\alpha|=-4,~~r\leq 1,~~l\leq-3$.\\

\indent By Theorem \ref{thm2}, we can compute the interior of $\widetilde{{\rm Wres}}[\pi^+(\nabla_X^{S(TM)}\nabla_Y^{S(TM)}D^{-1})\circ\pi^+(D^{-3})]$,\\
 we get
\begin{align}
\label{a16}
{\rm Wres}[\nabla_X^{S(TM)}\nabla_Y^{S(TM)}D^{-4}]
&=\frac{4\pi^2}{3}\int_{M}[Ric(X,Y)-\frac{1}{2}sg(X,Y)]d{\rm Vol_{M}}+\pi^2\int_{M}sg(X,Y)d{\rm Vol_{M}}.
\end{align}
\indent Now we  need to compute $\int_{\partial M} \Phi$. When $n=4$, then ${\rm tr}_{S(TM)}[{\rm \texttt{id}}]={\rm dim}(\wedge^*(\mathbb{R}^2))=4$, the sum is taken over $
r+l-k-j-|\alpha|=-3,~~r\leq 0,~~l\leq-2,$ then we have the following five cases:
~\\

\noindent  {\bf case a)~I)}~$r=1,~l=-2,~k=j=0,~|\alpha|=1$.\\
\noindent By (\ref{b15}), we get
\begin{equation}
\label{b24}
\Phi_1=-\int_{|\xi'|=1}\int^{+\infty}_{-\infty}\sum_{|\alpha|=1}
 {\rm tr}[\partial^\alpha_{\xi'}\pi^+_{\xi_n}\sigma_{1}(\nabla_X^{S(TM)}\nabla_Y^{S(TM)}D^{-1})\times
 \partial^\alpha_{x'}\partial_{\xi_n}\sigma_{-3}(D^{-3})](x_0)d\xi_n\sigma(\xi')dx'.
\end{equation}
By Lemma 2.2 in \cite{Wa3}, for $i<n$, then
\begin{equation}
\label{b25}
\partial_{x_i}\sigma_{-3}({D}^{-3})(x_0)=
\partial_{x_i}(\sqrt{-1}c(\xi)|\xi|^{-4})(x_0)=
\sqrt{-1}\frac{\partial_{x_i}c(\xi)}{|\xi|^4}(x_0)
+\sqrt{-1}\frac{c(\xi)\partial_{x_i}(|\xi|^4)}{|\xi|^8}(x_0)
=0,
\end{equation}
\noindent so $\Phi_1=0$.\\

 \noindent  {\bf case a)~II)}~$r=1,~l=-3,~k=|\alpha|=0,~j=1$.\\
\noindent By (\ref{b15}), we get
\begin{equation}
\label{b26}
\Phi_2=-\frac{1}{2}\int_{|\xi'|=1}\int^{+\infty}_{-\infty} {\rm
trace} [\partial_{x_n}\pi^+_{\xi_n}\sigma_{1}(\nabla_X^{S(TM)}\nabla_Y^{S(TM)}D^{-1})\times
\partial_{\xi_n}^2\sigma_{-3}(D^{-3})](x_0)d\xi_n\sigma(\xi')dx'.
\end{equation}

\noindent By Lemma \ref{lem3} and Lemma\ref{lem356}, we have\\
\begin{eqnarray}\label{b237}
\partial_{\xi_n}^2\sigma_{-3}((D^{-3}))(x_0)=\partial_{\xi_n}^2(c(\xi)|\xi|^{-4})(x_0)
=\sqrt{-1}\frac{(20\xi_n^2-4)c(\xi')+12(\xi^3-\xi)c(\mathrm{d}x_n)}{(1+\xi_n^2)^4},
\end{eqnarray}
and
\begin{align}\label{b27}
\partial_{x_n}\sigma_{1}(\nabla_X^{S(TM)}\nabla_Y^{S(TM)}D^{-1})&=\partial_{x_n}(-\sqrt{-1}\Sigma_{j,l=1}^nX_jY_l\xi_j\xi_lc(\xi)|\xi|^{-2})\nonumber\\
&=\Sigma_{j,l=1}^nX_jY_l\xi_j\xi_l \left[ \frac{\partial_{x_n}c(\xi')}{1+\xi_n^2}+\frac{c(\xi)h'(0)|\xi'|^2}{(1+\xi_n^2)^2}\right].
\end{align}

Then, we have
\begin{align}\label{b28}
\pi^+_{\xi_n}\partial_{x_n}\sigma_{1}(\nabla_X^{S(TM)}\nabla_Y^{S(TM)}D^{-1})
&=\partial_{x_n}\pi^+_{\xi_n}\sigma_{1}(\nabla_X^{S(TM)}\nabla_Y^{S(TM)}D^{-1})\nonumber\\
&=\sqrt{-1}\Sigma_{j,l=1}^{n-1}X_jY_l\xi_j\xi_lh'(0)|\xi'|^2\left[\frac{ic(\xi')}{4(\xi_n-i)}+\frac{c(\xi')+ic(\mathrm{d}x_n)}{4(\xi_n-i)^2}\right]\nonumber\\
&-\Sigma_{j,l=1}^{n-1}X_jY_l\xi_j\xi_l\frac{\partial_{x_n}c(\xi')}{2(\xi_n-i)}\nonumber\\
&-\sqrt{-1}X_nY_n\left\{\frac{\partial_{x_n}c(\xi')}{2(\xi_n-i)}
 +h'(0)|\xi'|\left[-\frac{2ic(\xi')-3c(\mathrm{d}x_n)}{4(\xi_n-i)}\right.\right.\nonumber\\
&+\left.\left.\frac{[c(\xi')+ic(\mathrm{d}x_n)][i(\xi_n-i)+1]}{4(\xi_n-i)^2}\right]\right\}\nonumber\\
&-\Sigma_{j=1}^{n-1}X_jY_n\xi_j\left[\sqrt{-1}\frac{\partial_{x_n}c(\xi')}{2(\xi_n-i)}
-\frac{\sqrt{-1}h'(0)|\xi'|[c(\xi')+2ic(\mathrm{d}x_n)]}{4(\xi_n-i)}\right.\nonumber\\
&\left.-\frac{[ic(\xi')-c(\mathrm{d}x_n)][i(\xi_n-i)+1]}{(\xi_n-i)^2}\right]\nonumber\\
&-\Sigma_{l=1}^{n-1}X_nY_l\xi_l \left[\sqrt{-1}\frac{\partial_{x_n}c(\xi')}{2(\xi_n-i)}
-\frac{\sqrt{-1}h'(0)|\xi'|[c(\xi')+2ic(\mathrm{d}x_n)]}{4(\xi_n-i)}\right.\nonumber\\
&\left.-\frac{[ic(\xi')-c(\mathrm{d}x_n)][i(\xi_n-i)+1]}{(\xi_n-i)^2}\right] .\nonumber\\
\end{align}

We note that $i<n,~\int_{|\xi'|=1}\xi_{i_{1}}\xi_{i_{2}}\cdots\xi_{i_{2d+1}}\sigma(\xi')=0$,
so we omit some items that have no contribution for computing {\bf case a)~II)}.\\
Then there is the following formula
\begin{align}\label{33}
&{\rm
tr} [\partial_{x_n}\pi^+_{\xi_n}\sigma_{1}(\nabla_X^{S(TM)}\nabla_Y^{S(TM)}D^{-1})\times
\partial_{\xi_n}^2\sigma_{-3}(D^{-3})](x_0)\nonumber\\
&=\Sigma_{j,l=1}^{n-1}X_jY_l\xi_j\xi_lh'(0)\left[8i\frac{5\xi_n^2-1}{(\xi_n-i)^5(\xi_n+i)^4}
 +\frac{4(5\xi_n^2-1)+12i(\xi_n^3-\xi_n)}{(\xi_n-i)^6(\xi_n+i)^4}\right]\nonumber\\
&+X_nY_nh'(0)\left[\frac{(4i-4)(\xi^2-1)+48(\xi_n^3-\xi_n)}{(\xi_n-i)^5(\xi_n+i)^4}
 -\frac{4(5\xi_n^2-1)+12i(\xi_n^3-\xi_n)}{(\xi_n-i)^6(\xi_n+i)^4}\right]\nonumber\\
&+8\Sigma_{j=1}^{n-1}X_jY_n\xi_j\left[\frac{(6-3ih'(0))(\xi_n^3-\xi_n)-2i(5\xi_n^2-1)}{(\xi_n-i)^5(\xi_n+i)^4}
 +\frac{2(5\xi_n^2-1)+6i(\xi_n^3-\xi_n)}{(\xi_n-i)^6(\xi_n+i)^4}\right]\nonumber\\
&+8\Sigma_{l=1}^{n-1}X_nY_l\xi_l\left[\frac{(6-3ih'(0))(\xi_n^3-\xi_n)-2i(5\xi_n^2-1)}{(\xi_n-i)^5(\xi_n+i)^4}
 +\frac{2(5\xi_n^2-1)+6i(\xi_n^3-\xi_n)}{(\xi_n-i)^6(\xi_n+i)^4}\right].\nonumber\\
\end{align}

Therefore, we get
\begin{align}\label{35}
\Phi_2&=\frac{1}{2}\int_{|\xi'|=1}\int^{+\infty}_{-\infty}\bigg\{
\Sigma_{j,l=1}^{n-1}X_jY_l\xi_j\xi_lh'(0)\left[8i\frac{5\xi_n^2-1}{(\xi_n-i)^5(\xi_n+i)^4}
 +\frac{4(5\xi_n^2-1)+12i(\xi_n^3-\xi_n)}{(\xi_n-i)^6(\xi_n+i)^4}\right]\nonumber\\
&+X_nY_nh'(0)\left[\frac{(4i-4)(\xi^2-1)+48(\xi_n^3-\xi_n)}{(\xi_n-i)^5(\xi_n+i)^4}
 -\frac{4(5\xi_n^2-1)+12i(\xi_n^3-\xi_n)}{(\xi_n-i)^6(\xi_n+i)^4}\right]
 \bigg\}d\xi_n\sigma(\xi')dx'\nonumber\\
 &=\Sigma_{j,l=1}^{n-1}X_jY_lh'(0)\Omega_3\int_{\Gamma^{+}}\left[8i\frac{5\xi_n^2-1}{(\xi_n-i)^5(\xi_n+i)^4}
 +\frac{4(5\xi_n^2-1)+12i(\xi_n^3-\xi_n)}{(\xi_n-i)^6(\xi_n+i)^4}\right]\xi_j\xi_ld\xi_{n}dx'\nonumber\\
 &+X_nY_nh'(0)\Omega_3\int_{\Gamma^{+}}\left[\frac{(4i-4)(\xi^2-1)+48(\xi_n^3-\xi_n)}{(\xi_n-i)^5(\xi_n+i)^4}
 -\frac{4(5\xi_n^2-1)+12i(\xi_n^3-\xi_n)}{(\xi_n-i)^6(\xi_n+i)^4}\right]d\xi_{n}dx'\nonumber\\
&=\Sigma_{j=1}^{n-1}X_jY_j\frac{4\pi}{3}h'(0)\Omega_3
\left\{ \frac{2\pi i}{4!}
\left[\frac{8i(5\xi_n^2-1)}{(\xi_n+i)^4}\right]^{(4)}\bigg|_{\xi_n=i}dx'
+\frac{2\pi i}{5!}
\left[\frac{4(5\xi_n^2-1)}{(\xi_n+i)^4}\right]^{(5)}\bigg|_{\xi_n=i}dx' \right.\nonumber\\
 &\left.+\frac{2\pi i}{5!}
\left[\frac{\xi_n^3-\xi_n}{(\xi_n+i)^4}\right]^{(5)}\bigg|_{\xi_n=i}dx'
\right\}\nonumber\\
 &-X_nY_nh'(0)\Omega_3\left\{ \frac{2\pi i}{4!}
 \left[\frac{(4i-4)(\xi_n^2-1)}{(\xi_n+i)^4}\right]^{(4)}\bigg|_{\xi_n=i}dx'
 +\frac{2\pi i}{4!}
 \left[\frac{48(\xi_n^3-\xi_n)}{(\xi_n+i)^4}\right]^{(4)}\bigg|_{\xi_n=i}dx'\right.\nonumber\\
 &\left.-\frac{2\pi i}{5!}
 \left[\frac{4(5\xi_n^2-1)}{(\xi_n+i)^4}\right]^{(5)}\bigg|_{\xi_n=i}dx'
 -\frac{2\pi i}{5!}
 \left[\frac{12i(\xi_n^3-\xi_n)}{(\xi_n+i)^4}\right]^{(5)}\bigg|_{\xi_n=i}dx'
\right\}\nonumber\\
&=-\left[\frac{592}{3}\pi\Sigma_{j=1}^{n-1}X_jY_j
+\left(\frac{461}{4}+\frac{23}{4}i\right)X_nY_n\right]h'(0)\pi\Omega_3dx',
\end{align}
where ${\rm \Omega_{3}}$ is the canonical volume of $S^{3}.$\\

\noindent  {\bf case a)~III)}~$r=1,~l=-3,~j=|\alpha|=0,~k=1$.\\
\noindent By (\ref{b15}), we get
\begin{align}\label{36}
\Phi_3&=-\frac{1}{2}\int_{|\xi'|=1}\int^{+\infty}_{-\infty}
{\rm trace} [\partial_{\xi_n}\pi^+_{\xi_n}\sigma_{1}(\nabla_X^{S(TM)}\nabla_Y^{S(TM)}D^{-1})\times
\partial_{\xi_n}\partial_{x_n}\sigma_{-3}(D^{-3})](x_0)d\xi_n\sigma(\xi')dx'\nonumber\\
&=\frac{1}{2}\int_{|\xi'|=1}\int^{+\infty}_{-\infty}
{\rm trace} [\partial_{\xi_n}^2\pi^+_{\xi_n}\sigma_{1}(\nabla_X^{S(TM)}\nabla_Y^{S(TM)}D^{-1})\times
\partial_{x_n}\sigma_{-3}(D^{-3})](x_0)d\xi_n\sigma(\xi')dx'.
\end{align}
\noindent By Lemma \ref{lem356}, we have
\begin{eqnarray}\label{37}
\partial_{x_n}\sigma_{-3}(D^{-3})(x_0)|_{|\xi'|=1}
=\frac{\sqrt{-1}\partial_{x_n}[c(\xi')]}{(1+\xi_n^2)^4}-\frac{2\sqrt{-1}h'(0)c(\xi)|\xi'|^2_{g^{\partial M}}}{(1+\xi_n^2)^6},
\end{eqnarray}
moreover,
\begin{align}\label{38}
\pi^+_{\xi_n}\sigma_{1}(\nabla_X^{S(TM)}\nabla_Y^{S(TM)}D^{-1})
&=-\frac{c(\xi')+ic(\mathrm{d}x_n)}{2(\xi_n-i)}\Sigma_{j,l=1}^{n-1}X_jY_l\xi_j\xi_l
-\frac{c(\xi')+ic(\mathrm{d}x_n)}{2(\xi_n-i)}X_nY_n\nonumber\\
&-\frac{ic(\xi')-c(\mathrm{d}x_n)}{2(\xi_n-i)}\Sigma_{j=1}^{n-1}X_jY_n\xi_j
-\frac{ic(\xi')-c(\mathrm{d}x_n)}{2(\xi_n-i)}\Sigma_{l=1}^{n-1}X_nY_l\xi_l.
\end{align}

Then, we have
\begin{align}\label{mmmmm}
\partial_{\xi_n}^2\pi^+_{\xi_n}\sigma_{1}(\nabla_X^{S(TM)}\nabla_Y^{S(TM)}D^{-1})
=-\frac{c(\xi')+ic(\mathrm{d}x_n)}{(\xi_n-i)^3}\Sigma_{j,l=1}^{n-1}X_jY_l\xi_j\xi_l
-\frac{c(\xi')+ic(\mathrm{d}x_n)}{(\xi_n-i)^3}X_nY_n.
\end{align}

We note that $i<n,~\int_{|\xi'|=1}\xi_{i_{1}}\xi_{i_{2}}\cdots\xi_{i_{2d+1}}\sigma(\xi')=0$,
so we omit some items that have no contribution for computing {\bf case a)~III)}\\

\begin{align}\label{39}
&{\rm tr} [\partial_{\xi_n}\pi^+_{\xi_n}\sigma_{1}(\nabla_X^{S(TM)}\nabla_Y^{S(TM)}D^{-1})\times
\partial_{\xi_n}\partial_{x_n}\sigma_{-3}(D^{-3})](x_0)\nonumber\\
&=-2\frac{h'(0)}{(\xi_n-i)^5(\xi_n+i)^2}\Sigma_{j,l=1}^{n-1}X_jY_l\xi_j\xi_l
-2\frac{h'(0)}{(\xi_n-i)^5(\xi_n+i)^2}X_nY_n.\nonumber\\
\end{align}

Therefore, we get
\begin{align}\label{41}
\Phi_3&=\frac{1}{2}\int_{|\xi'|=1}\int^{+\infty}_{-\infty}
\bigg(-2\frac{h'(0)}{(\xi_n-i)^5(\xi_n+i)^2}\Sigma_{j,l=1}^{n-1}X_jY_l\xi_j\xi_l
-2\frac{h'(0)}{(\xi_n-i)^5(\xi_n+i)^2}X_nY_n\bigg)d\xi_n\sigma(\xi')dx'\nonumber\\
&=-\Sigma_{j,l=1}^{n-1}X_jY_lh'(0)\Omega_3\int_{\Gamma^{+}}\frac{1}{(\xi_n-i)^5(\xi_n+i)^2}\xi_j\xi_ld\xi_{n}dx'
-X_nY_nh'(0)\Omega_3\int_{\Gamma^{+}}\frac{1}{(\xi_n-i)^5(\xi_n+i)^2}d\xi_{n}dx'\nonumber\\
&=-\Sigma_{j,l=1}^{n-1}X_jY_lh'(0)\Omega_3\frac{2\pi i}{4!}
\left[\frac{1}{(\xi_n+i)^2}\right]^{(4)}\bigg|_{\xi_n=i}dx'
+2X_nY_nh'(0)\Omega_3\frac{2\pi i}{4!}
\left[\frac{1}{(\xi_n+i)^2}\right]^{(4)}\bigg|_{\xi_n=i}dx'\nonumber\\
&=\left(\frac{5\pi i}{6}\Sigma_{j=1}^{n-1}X_jY_j+\frac{5i}{8}X_nY_n\right)h'(0)\pi\Omega_3dx'.
\end{align}

\noindent  {\bf case b)}~$r=0,~l=-3,~k=j=|\alpha|=0$.\\
\noindent By (\ref{b15}), we get
\begin{align}\label{42}
\Phi_4&=-i\int_{|\xi'|=1}\int^{+\infty}_{-\infty}{\rm trace} [\pi^+_{\xi_n}\sigma_{0}(\nabla_X^{S(TM)}\nabla_Y^{S(TM)}D^{-1})\times
\partial_{\xi_n}\sigma_{-3}(D^{-3})](x_0)d\xi_n\sigma(\xi')dx'\nonumber\\
&=i\int_{|\xi'|=1}\int^{+\infty}_{-\infty}{\rm trace} [\partial_{\xi_n}\pi^+_{\xi_n}\sigma_{0}(\nabla_X^{S(TM)}\nabla_Y^{S(TM)}D^{-1})\times
\sigma_{-3}(D^{-3})](x_0)d\xi_n\sigma(\xi')dx'.
\end{align}

 By Lemma \ref{lem356}, we have
\begin{align}\label{43}
\partial_{\xi_n}\sigma_{-3}(D^{-3})(x_0)|_{|\xi'|=1}=\frac{ic(\mathrm{d}x_n)}{(1+\xi_n^2)^2}-\frac{4\sqrt{-1}\xi_nc(\xi)}{(1+\xi_n^2)^3}.
\end{align}
By Lemma \ref{lem3}, we have
\begin{align}
\sigma_{0}(\nabla_X^{S(TM)}\nabla_Y^{S(TM)}D^{-1})&=-\Sigma_{j,l=1}^{n}X_jY_l\xi_j\xi_l\left[\frac{c(\xi)\sigma_0(D)c(\xi)}{|\xi|^4}
+\frac{c(\xi)}{|\xi|^6}\Sigma_jc(\mathrm{d}x_j)(\partial_{x_j}[c(\xi)]|\xi|^2-c(\xi)\partial_{x_j}(|\xi|^2))\right]\nonumber\\
&-i\Sigma_{j,l=1}^nX_j\frac{\partial_{Y_l}}{\partial_{x_j}}\xi_j\frac{c(\xi)}{|\xi|^2}
-\Sigma_jA(Y)X_j\xi_j\frac{c(\xi)}{|\xi|^2}-\Sigma_lA(Y)Y_l\xi_l\frac{c(\xi)}{|\xi|^2},
\end{align}
moreover,
\begin{align}\label{45}
\partial_{\xi_n}\pi^+_{\xi_n}\sigma_{0}(\nabla_X^{S(TM)}\nabla_Y^{S(TM)}D^{-1})
&=-\Sigma_{j,l=1}^{n-1}X_jY_l\xi_j\xi_l[B_1-B_2]\nonumber\\
&-X_nY_n\left\{\sigma_0\left[\frac{[c(\xi')+ic(\mathrm{d}x_n)]c(\mathrm{d}x_n)}{4(\xi_n-i)}
  +\frac{[c(\xi')+ic(\mathrm{d}x_n)]^2}{4(\xi_n-i)^2}\right]\right.\nonumber\\
&+c(\mathrm{d}x_n)\partial_{x_n}[c(\xi')]\left[\frac{i[c(\xi')+ic(\mathrm{d}x_n)]c(\mathrm{d}x_n)}{-4(\xi_n-i)}
  +\frac{[c(\xi')+ic(\mathrm{d}x_n)]}{4(\xi_n-i)^2}\right] \nonumber\\
&+h'(0)|\xi'|^2\left[\frac{i[c(\xi')+ic(\mathrm{d}x_n)]^2-2c(\mathrm{d}x_n)[3c(\xi')+4ic(\mathrm{d}x_n)]}{16(\xi_n-i)}\right.\nonumber\\ &\left.\left.+\frac{[c(\xi')+ic(\mathrm{d}x_n)][c(\xi')+5ic(\mathrm{d}x_n)]}{16(\xi_n-i)^2}
-\frac{i[c(\xi')+ic(\mathrm{d}x_n)]^2}{8(\xi_n-i)^3}\right]  \right\}\nonumber\\
&-\Sigma_{j=1}^{n-1}X_jY_n\xi_j\left\{\sigma_0\left[\frac{[-ic(\xi')+ic(\mathrm{d}x_n)]c(\mathrm{d}x_n)}{2(\xi_n-i)}
  -\frac{i[c(\xi')+ic(\mathrm{d}x_n)]^2}{4(\xi_n-i)^2}\right]\right.\nonumber\\
&+c(\mathrm{d}x_n)\partial_{x_n}[c(\xi')]\left[\frac{-ic(\mathrm{d}x_n)}{4(\xi_n-i)}
  -\frac{i[c(\xi')+ic(\mathrm{d}x_n)]}{4(\xi_n-i)^2}\right] \nonumber\\
&-c(\mathrm{d}x_n)h'(0)|\xi'|^2\left[\frac{ic(\mathrm{d}x_n)c(\xi')}{8(\xi_n-i)}
 +\frac{[c(\xi')+ic(\mathrm{d}x_n)][c(\xi')+3ic(\mathrm{d}x_n)]}{8(\xi_n-i)^2}\right.\nonumber\\
&\left.\left.-\frac{[c(\xi')+ic(\mathrm{d}x_n)]^2}{8(\xi_n-i)^3}\right]  \right\}\nonumber\\
&-\Sigma_{l=1}^{n-1}X_nY_l\xi_l\left\{\sigma_0\left[\frac{[-ic(\xi')+ic(\mathrm{d}x_n)]c(\mathrm{d}x_n)}{2(\xi_n-i)}
  -\frac{i[c(\xi')+ic(\mathrm{d}x_n)]^2}{4(\xi_n-i)^2}\right]\right.\nonumber\\
&+c(\mathrm{d}x_n)\partial_{x_n}[c(\xi')]\left[\frac{-ic(\mathrm{d}x_n)}{4(\xi_n-i)}
  -\frac{i[c(\xi')+ic(\mathrm{d}x_n)]}{4(\xi_n-i)^2}\right] \nonumber\\
&-c(\mathrm{d}x_n)h'(0)|\xi'|^2\left[\frac{ic(\mathrm{d}x_n)c(\xi')}{8(\xi_n-i)}
 +\frac{[c(\xi')+ic(\mathrm{d}x_n)][c(\xi')+3ic(\mathrm{d}x_n)]}{8(\xi_n-i)^2}\right.\nonumber\\
&\left.\left.-\frac{[c(\xi')+ic(\mathrm{d}x_n)]^2}{8(\xi_n-i)^3}\right]  \right\}\nonumber\\
&-\Sigma_{j,l=1}^nX_j\frac{\partial_{Y_l}}{\partial_{x_j}}\xi_j\frac{c(\xi')+ic(\mathrm{d}x_n)}{2(\xi_n-i)^2}
-i\Sigma_jA(Y)X_j\xi_j\frac{c(\xi')+ic(\mathrm{d}x_n)}{2(\xi_n-i)^2}\nonumber\\
&-i\Sigma_lA(Y)Y_l\xi_l\frac{c(\xi')+ic(\mathrm{d}x_n)}{2(\xi_n-i)^2}.
\end{align}

Similarly to \cite{Wa3}, we have

\begin{align}
B_1=-\frac{A_1}{4(\xi_n-i)}-\frac{A_2}{4(\xi_n-i)^2}
\end{align}
where
\begin{align}
A_1&=ic(\xi')p_0c(\xi')+ic(\mathrm{d}x_n)p_0c(\mathrm{d}x_n)+ic(\xi')c(\mathrm{d}x_n)\partial_{x_n}[c(\xi')];\nonumber\\
A_2&=[c(\xi')+ic(\mathrm{d}x_n)]p_0[c(\xi')+ic(\mathrm{d}x_n)]+c(\xi')c(\mathrm{d}x_n)\partial_{x_n}[c(\xi')]-i\partial_{x_n}[c(\xi')];\\
B_2&=\frac{h'(0)}{2}\left[\frac{c(\mathrm{d}x_n)}{4i(\xi_n-i)}+\frac{c(\mathrm{d}x_n)-ic(\xi')}{8(\xi_n-i)^2}+\frac{3\xi_n-7i}{8(\xi_n-i)^3}[ic(\xi')-c(\mathrm{d}x_n)]\right].
\end{align}

We note that $i<n,~\int_{|\xi'|=1}\xi_{i_{1}}\xi_{i_{2}}\cdots\xi_{i_{2d+1}}\sigma(\xi')=0$,
so we omit some items that have no contribution for computing {\bf case b)}.\\
Then, we have
\begin{align}\label{39}
&{\rm tr}[\partial_{\xi_n}\pi^+_{\xi_n}\sigma_{0}(\nabla_X^{S(TM)}\nabla_Y^{S(TM)}D^{-1})\times
\sigma_{-3}(D^{-3})](x_0)\nonumber\\
&=\left[\frac{-2ic_0}{(1+\xi_n^2)^2}+h'(0)\frac{\xi_n^2-i\xi_n-2}{2(\xi_n-i)(1+\xi_n^2)^2}
-\frac{\sqrt{-1}}{2}h'(0)\frac{-i\xi_n^2-\xi_n+4i}{4(\xi_n-i)^3(\xi_n+i)^2}tr[id]\right]\Sigma_{j,l=1}^{n-1}X_jY_l\xi_j\xi_l\nonumber\\
&+\left[\frac{4-4i\xi_n-3(\xi_n^2+1)}{(\xi_n-i)^4(\xi_n+i)^3}+\frac{4i+4\xi_n-3i(\xi_n^2+1)}{2(\xi_n-i)^5(\xi_n+i)^3}\right]h'(0)X_nY_n.\nonumber\\
\end{align}

Therefore, we get
\begin{align}\label{41}
\Phi_4&=i\int_{|\xi'|=1}\int^{+\infty}_{-\infty}
\left[\frac{-2ic_0}{(1+\xi_n^2)^2}+h'(0)\frac{\xi_n^2-i\xi_n-2}{2(\xi_n-i)(1+\xi_n^2)^2}\right.\nonumber\\
&-\left.\frac{\sqrt{-1}}{2}h'(0)\frac{-i\xi_n^2-\xi_n+4i}{(\xi_n-i)^3(\xi_n+i)^2}\right]\Sigma_{j,l=1}^{n-1}X_jY_l\xi_j\xi_l\nonumber\\
&+\left[\frac{4-4i\xi_n-3(\xi_n^2+1)}{(\xi_n-i)^4(\xi_n+i)^3}+\frac{4i+4\xi_n-3i(\xi_n^2+1)}{2(\xi_n-i)^5(\xi_n+i)^3}\right]h'(0)X_nY_nd\xi_n\sigma(\xi')dx'\nonumber\\
&=\Sigma_{j,l=1}^{n-1}X_jY_lh'(0)\Omega_3\int_{\Gamma^{+}}\left[\frac{-2ic_0}{(1+\xi_n^2)^2}+h'(0)\frac{\xi_n^2-i\xi_n-2}{2(\xi_n-i)(1+\xi_n^2)^2}\right.\nonumber\\
&-\left.\frac{\sqrt{-1}}{2}h'(0)\frac{-i\xi_n^2-\xi_n+4i}{(\xi_n-i)^3(\xi_n+i)^2}\right]\xi_j\xi_ld\xi_{n}dx'\nonumber\\
&+X_nY_nh'(0)\Omega_3\int_{\Gamma^{+}}\left[\frac{4-4i\xi_n-3(\xi_n^2+1)}{(\xi_n-i)^4(\xi_n+i)^3}+\frac{4i+4\xi_n-3i(\xi_n^2+1)}{2(\xi_n-i)^5(\xi_n+i)^3}\right]d\xi_{n}dx'\nonumber\\
&=\frac{55}{3}\Sigma_{j,l=1}^{n-1}X_jY_lh'(0)\pi^2\Omega_3dx'
+X_nY_nh'(0)\Omega_3i\left\{\frac{2\pi i}{3!}\left[\frac{1}{(\xi_n+i)^3}\right]^{(3)}\bigg|_{\xi_n=i}dx'\right.\nonumber\\
&-\frac{4i\times2\pi i}{3!}\left[\frac{\xi_n}{(\xi_n+i)^3}\right]^{(3)}\bigg|_{\xi_n=i}dx'
-\frac{3\times2\pi i}{2!}\left[\frac{1}{(\xi_n+i)^2}\right]^{(2)}\bigg|_{\xi_n=i}dx'\nonumber\\
&+\frac{4i\times2\pi i}{2\times4!}\left[\frac{1}{(\xi_n+i)^3}\right]^{(4)}\bigg|_{\xi_n=i}dx'
+\frac{4\times2\pi i}{2\times4!}\left[\frac{\xi_n}{(\xi_n+i)^3}\right]^{(4)}\bigg|_{\xi_n=i}dx'\nonumber\\
&-\left.\frac{3i\times2\pi i}{2\times3!}\left[\frac{1}{(\xi_n+i)^2}\right]^{(3)}\bigg|_{\xi_n=i}dx'\right\}\nonumber\\
&=\left(\frac{55\pi}{3}\Sigma_{j=1}^{n-1}X_jY_j-\frac{3}{8}X_nY_n\right)h'(0)\pi\Omega_3dx'.
\end{align}

\noindent {\bf  case c)}~$r=1,~\ell=-4,~k=j=|\alpha|=0$.\\
By (\ref{b15}), we get
\begin{align}\label{61}
\Phi_5&=-\int_{|\xi'|=1}\int^{+\infty}_{-\infty}{\rm trace} [\pi^+_{\xi_n}\sigma_{1}(\nabla_X^{S(TM)}\nabla_Y^{S(TM)}D^{-1})\times
\partial_{\xi_n}\sigma_{-4}(D^{-3})](x_0)d\xi_n\sigma(\xi')dx'\nonumber\\
&=\int_{|\xi'|=1}\int^{+\infty}_{-\infty}{\rm trace} [\partial_{\xi_n}\pi^+_{\xi_n}\sigma_{1}(\nabla_X^{S(TM)}\nabla_Y^{S(TM)}D^{-1})\times
\sigma_{-4}(D^{-3})](x_0)d\xi_n\sigma(\xi')dx'.
\end{align}

By Lemma\ref{lem3} and \ref{lem356}, we have
\begin{align}\label{62}
\sigma_{-4}(D^{-3})=&\frac{1}{(\xi_n^2+1)^4}
\left[\left(\frac{11}{2}\xi_n(1+\xi_n^2)+8i\xi_n\right)h'(0)c(\xi')\right.\nonumber\\
&+\left[-2i+6i\xi_n^2-\frac{7}{4}(1+\xi_n^2)
 +\frac{15}{4}\xi_n^2(1+\xi^2_n)\right]h'(0)c(\mathrm{d}x_n) \nonumber\\
&\left.-3i\xi_n(1+\xi^2_n)\partial_{x_n}c(\xi')
 +i(1+\xi^2_n)c(\xi')c(\mathrm{d}x_n)\partial_{x_n}c(\xi')\right],
\end{align}
and
\begin{align}\label{621}
\partial_{\xi_n}\pi^+_{\xi_n}\sigma_{1}(\nabla_X^{S(TM)}\nabla_Y^{S(TM)}D^{-1})
&=\frac{c(\xi')+ic(\mathrm{d}x_n)}{2(\xi_n-i)^2}\Sigma_{j,l=1}^{n-1}X_jY_l\xi_j\xi_l
-\frac{c(\xi')+ic(\mathrm{d}x_n)}{2(\xi_n-i)^2}X_nY_n \nonumber\\
&+\frac{ic(\xi')-c(\mathrm{d}x_n)}{2(\xi_n-i)^2}\Sigma_{j=1}^{n}X_jY_n\xi_j
+\frac{ic(\xi')-c(\mathrm{d}x_n)}{2(\xi_n-i)^2}\Sigma_{l=1}^{n}X_nY_l\xi_l.
\end{align}

We note that $i<n,~\int_{|\xi'|=1}\xi_{i_{1}}\xi_{i_{2}}\cdots\xi_{i_{2d+1}}\sigma(\xi')=0$,
so we omit some items that have no contribution for computing {\bf case c)}.\\
We have
\begin{align}\label{63}
{\rm tr}[c(\xi')c(\xi')c(\mathrm{d}x_n)\partial_{x_n}c(\xi')]=0 ;\nonumber\\
{\rm tr}[c(\mathrm{d}x_n)c(\xi')c(\mathrm{d}x_n)\partial_{x_n}c(\xi')]=-2h'(0)
\end{align}
Then we get
\begin{align}\label{71}
&{\rm tr}[\partial_{\xi_n}\pi^+_{\xi_n}\sigma_{-1}(\nabla_X^{S(TM)}\nabla_Y^{S(TM)}D^{-1})\times
\sigma_{-4}(D^{-3})](x_0)|_{|\xi'|=1}\nonumber\\
&=\frac{h'(0)(7+6i-(20-15i)\xi_n-(7-6i)\xi_n^2+15i\xi_n^3)}{(\xi_n-i)^5(\xi_n+i)^4}\sum_{j,l=1}^{n-1}X_jY_l\xi_j\xi_l\nonumber\\
&+\frac{(3i-11)\xi_n(1-\xi_n^2)-16i\xi_n+(13+\frac{7}{2}i)(1+\xi_n^2)-16-\frac{15}{2}\xi_n^2(1+\xi_n^2)}{(\xi_n-i)^2(\xi_n+i)^4}X_nY_n.\nonumber\\
\end{align}

So we have
\begin{align}\label{74}
\Phi_5&=\int_{|\xi'|=1}\int^{+\infty}_{-\infty}{\rm tr}[\partial_{\xi_n}\pi^+_{\xi_n}\sigma_{1}(\nabla_X^{S(TM)}\nabla_Y^{S(TM)}D^{-1})\times
\sigma_{-4}(D^{-3})](x_0)d\xi_n\sigma(\xi')dx'\nonumber\\
&=\int_{|\xi'|=1}\int^{+\infty}_{-\infty}\frac{h'(0)(7+6i-(20-15i)\xi_n-(7-6i)\xi_n^2+15i\xi_n^3)}{(\xi_n-i)^5(\xi_n+i)^4}\sum_{j,l=1}^{n-1}X_jY_l\xi_j\xi_l\nonumber\\
&+\frac{(3i-11)\xi_n(1-\xi_n^2)-16i\xi_n+(13+\frac{7}{2}i)(1+\xi_n^2)-16-\frac{15}{2}\xi_n^2(1+\xi_n^2)}{(\xi_n-i)^2(1+\xi_n^2)^4}X_nY_nd\xi_n\sigma(\xi')dx'\nonumber\\
&=\left(-\frac{35}{3}+\frac{50}{3}i\right)\sum_{j=1}^{n-1}X_jY_jh'(0)\pi^2\Omega_3dx' \nonumber\\
&+h'(0)\Omega_3\int_{\Gamma^{+}}\frac{(3i-11)\xi_n(1-\xi_n^2)-16i\xi_n+(13+\frac{7}{2}i)(1+\xi_n^2)-16-\frac{15}{2}\xi_n^2(1+\xi_n^2)}{(\xi_n-i)^6(\xi_n+i)^4}\sum_{j,l=1}^{n-1}X_jY_l\xi_j\xi_ld\xi_{n}dx'\nonumber\\
&=\left(-\frac{35}{3}+\frac{50}{3}i\right)\sum_{j=1}^{n-1}X_jY_jh'(0)\pi^2\Omega_3dx' \nonumber\\
&+h'(0)\Omega_3\frac{2\pi i}{5!}X_nY_n\left[\frac{(3i-11)\xi_n(1-\xi_n^2)-16i\xi_n+(13+\frac{7}{2}i)(1+\xi_n^2)-16-\frac{15}{2}\xi_n^2(1+\xi_n^2)}{(\xi_n+i)^4}\right]^{(5)}\bigg|_{\xi_n=i}dx'\nonumber\\
&=\left[\left(-\frac{35}{3}+\frac{50}{3}i\right)\sum_{j=1}^{n-1}X_jY_j\pi+\left(5-\frac{137}{32}i\right)X_nY_n\right]h'(0)\pi^2\Omega_3dx'.
\end{align}

Let $X=X^T+X_n\partial_n,~Y=Y^T+Y_n\partial_n,$ then we have $\sum_{j=1}^{n-1}X_jY_j=g(X^T,Y^T).$ Now $\Phi$ is the sum of the cases (a), (b) and (c). Therefore, we get
\begin{align}\label{795}
\Phi=\sum_{i=1}^5\Phi_i=\left[\left(-\frac{1411}{12}+\frac{27i}{32}\right)X_nY_n+\left(-\frac{572}{3}+\frac{35i}{2}\right)\pi g(X^T,Y^T)\right]h'(0)\pi\Omega_3dx'.
\end{align}
Therefore, by (\ref{b15})-(\ref{795}), we obtain following theorem

\begin{thm}\label{thmb1}
Let $M$ be a $4$-dimensional oriented
compact spin manifold with boundary $\partial M$ and the metric
$g^{M}$ be defined as (\ref{b1}), then we get the following equality:
\begin{align}
\label{b2773}
&\widetilde{{\rm Wres}}[\pi^+(\nabla_X^{S(TM)}\nabla_Y^{S(TM)}D^{-2})\circ\pi^+(D^{-2})]\nonumber\\
&=\frac{4\pi^2}{3}\int_{M}[Ric(X,Y)-\frac{1}{2}sg(X,Y)]d{\rm Vol_{M}}+\pi^2\int_{M}sg(X,Y)d{\rm Vol_{M}}\nonumber\\
&+\int_{\partial M}\left[\left(-\frac{1411}{12}+\frac{27i}{32}\right)X_nY_n+\left(-\frac{572}{3}+\frac{35i}{2}\right)\pi g(X^T,Y^T)\right]h'(0)\pi\Omega_3d{\rm Vol_{M}}.
\end{align}
\end{thm}

\section{The noncommutative residue for $3$-dimensional manifolds with boundary}
\label{section:4}
In this section, we compute the noncommutative residue $\widetilde{{\rm Wres}}[\pi^+(\nabla_X^{S(TM)}\nabla_Y^{S(TM)}D^{-1})\circ\pi^+(D^{-3})]$ on $3$-dimensional oriented compact manifolds with boundary. As in \cite{Wa1}, for an odd-dimensional manifolds with boundary, we have the formula
\begin{align}
\widetilde{{\rm Wres}}[\pi^+(\nabla_X^{S(TM)}\nabla_Y^{S(TM)}D^{-1})\circ\pi^+(D^{-3})]=\int_{\partial M}\Phi
\end{align}
 \indent Let $M$ be a 3-dimensional compact oriented manifold with boundary $\partial M$. We have
\begin{align}
{\rm tr}[c^2(\xi')]&=-2 ;
~{\rm tr}[c(\xi')c(\mathrm{d}x_n))]=0 ;  \nonumber\\
{\rm tr}[c(\mathrm{d}x_n))c(\mathrm{d}x_n))]&=-2 ;
~{\rm tr}[c(\xi')+ic(\mathrm{d}x_n)][c(\xi')+ic(\mathrm{d}x_n)]=-2-2i\xi_n.
\end{align}
By (\ref{b15}), satisfies $r+l-k-|\alpha|-j-1=-3$,
so we get ~$r=1,~\ell=-3,~k=j=|\alpha|=0$, then\\
\begin{align}
\Phi&=-\int_{|\xi'|=1}\int^{+\infty}_{-\infty}{\rm trace} [\pi^+_{\xi_n}\sigma_{1}(\nabla_X^{S(TM)}\nabla_Y^{S(TM)}D^{-1})\times
\partial_{\xi_n}\sigma_{-3}(D^{-3})](x_0)d\xi_n\sigma(\xi')dx' \nonumber\\
&=\int_{|\xi'|=1}\int^{+\infty}_{-\infty}{\rm trace} [\partial_{\xi_n}\pi^+_{\xi_n}\sigma_{1}(\nabla_X^{S(TM)}\nabla_Y^{S(TM)}D^{-1})\times
\sigma_{-3}(D^{-3})](x_0)d\xi_n\sigma(\xi')dx'.
\end{align}

We note that $i<n,~\int_{|\xi'|=1}\xi_{i_{1}}\xi_{i_{2}}\cdots\xi_{i_{2d+1}}\sigma(\xi')=0$,
so we omit some items that have no contribution for computing. By (\ref{621}) and (\ref{b22222}), we have

\begin{align}
&{\rm tr}[\partial_{\xi_n}\pi^+_{\xi_n}\sigma_{1}(\nabla_X^{S(TM)}\nabla_Y^{S(TM)}D^{-1})\times
\sigma_{-3}(D^{-3})]\nonumber\\
&=-\frac{i-\xi_n}{(\xi_n-i)^4(\xi_n+i)^2}\Sigma_{j,l=1}^{n-1}X_jY_l\xi_j\xi_l
-\frac{\sqrt{-1}\xi_n^2(2+2i\xi_n)}{2(\xi_n-i)^3(\xi_n+i)^2}X_nY_n.
\end{align}
Then there is the following formula
\begin{align}
\Phi&=-\int_{|\xi'|=1}\int^{+\infty}_{-\infty}{\rm trace} [\pi^+_{\xi_n}\sigma_{1}(\nabla_X^{S(TM)}\nabla_Y^{S(TM)}D^{-1})\times
\partial_{\xi_n}\sigma_{-3}(D^{-3})](x_0)d\xi_n\sigma(\xi')dx' \nonumber\\
&=\int_{|\xi'|=1}\int^{+\infty}_{-\infty}\frac{i-\xi_n}{(\xi_n-i)^4(\xi_n+i)^2}\Sigma_{j,l=1}^{n-1}X_jY_l\xi_j\xi_l
+\frac{\sqrt{-1}\xi_n^2(2+2i\xi_n)}{2(\xi_n-i)^3(\xi_n+i)^2}X_nY_n(x_0)d\xi_n\sigma(\xi')dx' \nonumber\\
&=\Omega_2\int_{\Gamma^{+}}-\frac{i-\xi_n}{(\xi_n-i)^4(\xi_n+i)^2}\Sigma_{j,l=1}^{n-1}X_jY_l\xi_j\xi_l
+\frac{\sqrt{-1}\xi_n^2(2+2i\xi_n)}{2(\xi_n-i)^3(\xi_n+i)^2}X_nY_n(x_0)d\xi_ndx' \nonumber\\
&=-\Omega_2\frac{2\pi i}{3!}\frac{4\pi}{3}\Sigma_{j,l=1}^{n-1}X_jY_l\left[\frac{i-\xi_n}{(\xi_n+i)^2}\right]^{(3)}\bigg|_{(\xi_n=i)}\mathrm{d}x'
+\Omega_2\frac{2\pi i}{3!}X_nY_n\left[\frac{\xi_n-i}{(\xi_n+i)^2}\right]^{(3)}\bigg|_{(\xi_n=i)}\mathrm{d}x'
\nonumber\\
&=\left(\frac{i}{3}\Sigma_{j=1}^{n-1}X_jY_j\pi+\frac{3i}{4}X_nY_n\right)\pi\Omega_2\mathrm{d}x'.
\end{align}
Therefore, we get the following theorem
\begin{thm}\label{thmb2}
Let $M$ be a $3$-dimensional oriented
compact spin manifold with boundary $\partial M$ and the metric
$g^{M}$ be defined as Sect.\ref{section:3}, then we get the following equality:
\begin{align}
\widetilde{{\rm Wres}}[\pi^+(\nabla_X^{S(TM)}\nabla_Y^{S(TM)}D^{-1})\circ\pi^+(D^{-3})]=
\left(\frac{i}{3}\Sigma_{j=1}^{n-1}X_jY_j\pi+\frac{3i}{4}X_nY_n\right)\pi\Omega_2\mathrm{d}x'.
\end{align}
\end{thm}


\section{The noncommutative residue for $6$-dimensional manifolds with boundary}
\label{section:5}
 In this section, we compute the noncommutative residue $\widetilde{{\rm Wres}}[\pi^+(\nabla_X^{S(TM)}\nabla_Y^{S(TM)}D^{-3})\circ\pi^+(D^{-3})]$ on $6$-dimensional oriented compact manifolds with boundary.  For more details, see Section 2 in \cite{Wa3}.\\
 \indent Let $M$ be a 6-dimensional compact oriented manifold with boundary $\partial M$ and the metric as same as Section \ref{section:3}.
 We can compute the noncommutative residue
\begin{align}
\label{c1}
&\widetilde{{\rm Wres}}[\pi^+(\nabla_X^{S(TM)}\nabla_Y^{S(TM)}D^{-3})\circ\pi^+(D^{-3})]\nonumber\\
&=\int_M\int_{|\xi|=1}{\rm
trace}_{\wedge^*T^*M\bigotimes\mathbb{C}}[\sigma_{-6}(\nabla_X^{S(TM)}\nabla_Y^{S(TM)}D^{-3}\circ D^{-3})]\sigma(\xi)dx+\int_{\partial M}\Phi,
\end{align}
where
\begin{align}
\label{c12}
\Phi &=\int_{|\xi'|=1}\int^{+\infty}_{-\infty}\sum^{\infty}_{j, k=0}\sum\frac{(-i)^{|\alpha|+j+k+1}}{\alpha!(j+k+1)!}
\times {\rm trace}_{\wedge^*T^*M\bigotimes\mathbb{C}}[\partial^j_{x_n}\partial^\alpha_{\xi'}\partial^k_{\xi_n}\sigma^+_{r}(\nabla_X^{S(TM)}\nabla_Y^{S(TM)}D^{-3})(x',0,\xi',\xi_n)
\nonumber\\
&\times\partial^\alpha_{x'}\partial^{j+1}_{\xi_n}\partial^k_{x_n}\sigma_{l}(D^{-3})(x',0,\xi',\xi_n)]d\xi_n\sigma(\xi')dx',
\end{align}
and the sum is taken over $r+l-k-j-|\alpha|=-5,~~r\leq -1,~~l\leq-3$.\\

\indent By Theorem \ref{thm2}, we can compute the interior of $\widetilde{{\rm Wres}}[\pi^+(\nabla_X^{S(TM)}\nabla_Y^{S(TM)}D^{-3})\circ\pi^+(D^{-3})]$,\\
 we get
\begin{align}
\label{c13}
{\rm Wres}[\nabla_X^{S(TM)}\nabla_Y^{S(TM)}D^{-6}]
&=\frac{4\pi^3}{3}\int_{M}[Ric(X,Y)-\frac{1}{2}sg(X,Y)]d{\rm Vol_{M}}+\pi^3\int_{M}sg(X,Y)d{\rm Vol_{M}}.
\end{align}
\indent Now we  need to compute $\int_{\partial M} \Phi$. When $n=6$, then ${\rm tr}_{S(TM)}[{\rm \texttt{id}}]={\rm dim}(\wedge^*(\mathbb{R}^3))=8$, the sum is taken over $
r+l-k-j-|\alpha|=-6,~~r\leq -1,~~l\leq-3,$ then we have the following five cases:
~\\

\noindent  {\bf case a)~I)}~$r=-1,~l=-3,~k=j=0,~|\alpha|=1$.\\
\noindent By (\ref{c12}), we get
\begin{equation}
\label{c1}
\Phi_1=-\int_{|\xi'|=1}\int^{+\infty}_{-\infty}\sum_{|\alpha|=1}
 {\rm tr}[\partial^\alpha_{\xi'}\pi^+_{\xi_n}\sigma_{-1}(\nabla_X^{S(TM)}\nabla_Y^{S(TM)}D^{-3})\times
 \partial^\alpha_{x'}\partial_{\xi_n}\sigma_{-3}(D^{-3})](x_0)d\xi_n\sigma(\xi')dx'.
\end{equation}
By Lemma 2.2 in \cite{Wa3}, for $i<n$, then
\begin{equation}
\label{c14}
\partial_{x_i}\sigma_{-3}({D}^{-3})(x_0)=
\partial_{x_i}(\sqrt{-1}c(\xi)|\xi|^{-4})(x_0)=
\sqrt{-1}\frac{\partial_{x_i}c(\xi)}{|\xi|^4}(x_0)
+\sqrt{-1}\frac{c(\xi)\partial_{x_i}(|\xi|^4)}{|\xi|^8}(x_0)
=0,
\end{equation}
\noindent so $\Phi_1=0$.\\

 \noindent  {\bf case a)~II)}~$r=-1,~l=-3,~k=|\alpha|=0,~j=1$.\\
\noindent By (\ref{c12}), we get
\begin{equation}
\label{c2}
\Phi_2=-\frac{1}{2}\int_{|\xi'|=1}\int^{+\infty}_{-\infty} {\rm
trace} [\partial_{x_n}\pi^+_{\xi_n}\sigma_{-1}(\nabla_X^{S(TM)}\nabla_Y^{S(TM)}D^{-3})\times
\partial_{\xi_n}^2\sigma_{-3}(D^{-3})](x_0)d\xi_n\sigma(\xi')dx'.
\end{equation}

\noindent By Lemma \ref{lem356}, we have\\
\begin{eqnarray}\label{c21}
\partial_{\xi_n}^2\sigma_{-3}((D^{-3}))(x_0)=\partial_{\xi_n}^2(c(\xi)|\xi|^{-4})(x_0)
=\sqrt{-1}\frac{(20\xi_n^2-4)c(\xi')+12(\xi^3-\xi)c(\mathrm{d}x_n)}{(1+\xi_n^2)^4}.
\end{eqnarray}

\begin{eqnarray}\label{c22}
\partial_{x_n}\sigma_{-1}(\nabla_X^{S(TM)}\nabla_Y^{S(TM)}D^{-3})
=\Sigma_{j,l=1}^nX_jY_l\xi_j\xi_l \left[ \frac{\partial_{x_n}c(\xi')}{1+\xi_n^2}^2-\frac{2c(\xi)h'(0)|\xi'|^2}{(1+\xi_n^2)^3}\right].
\end{eqnarray}

Then, we have
\begin{align}\label{c23}
\pi^+_{\xi_n}\partial_{x_n}\sigma_{-1}(\nabla_X^{S(TM)}\nabla_Y^{S(TM)}D^{-3})
&=\partial_{x_n}\pi^+_{\xi_n}\sigma_{-1}(\nabla_X^{S(TM)}\nabla_Y^{S(TM)}D^{-3})\nonumber\\
&=\sqrt{-1}\Sigma_{j,l=1}^{n-1}X_jY_l\xi_j\xi_lh'(0)|\xi'|^2\left[\frac{-3ic(\xi')}{4(\xi_n-i)}
-\frac{3c(\xi')+ic(\mathrm{d}x_n)}{4(\xi_n-i)^2}\right.\nonumber\\
&\left.+\frac{c(\xi')+ic(\mathrm{d}x_n)}{(\xi_n-i)^3}\right]
-\Sigma_{j,l=1}^{n-1}X_jY_l\xi_j\xi_l\frac{\partial_{x_n}c(\xi')[i(\xi_n-i)+1]}{4(\xi_n-i)^2}\nonumber\\
&-\sqrt{-1}X_nY_n\left\{\frac{-3i\partial_{x_n}c(\xi')}{4(\xi_n-i)}+\frac{\partial_{x_n}c(\xi')}{4(\xi_n-i)^2}
+\frac{c(\xi')+ic(\mathrm{d}x_n)}{2i(\xi_n-i)^3}\right.\nonumber\\
&\left.+h'(0)|\xi'|\left[-\frac{-5c(\xi')+9ic(\mathrm{d}x_n)}{4(\xi_n-i)}+\frac{3c(\xi')+3ic(\mathrm{d}x_n)}{4(\xi_n-i)^2}\right]
\right\}\nonumber\\
&-\Sigma_{j=1}^{n-1}X_jY_n\xi_j\left\{-\frac{i\partial_{x_n}c(\xi')}{4(\xi_n-i)^2}
+h'(0)|\xi'|\left[\frac{11ic(\mathrm{d}x_n)}{8(\xi_n-i)}\right.\right.\nonumber\\
&\left.\left.-\frac{ic(\xi')+(1+4i)c(\mathrm{d}x_n)}{16(\xi_n-i)^2}-\frac{ic(\xi')-c(\mathrm{d}x_n)}{2(\xi_n-i)^3}\right]\right\}\nonumber\\
&-\Sigma_{l=1}^{n-1}X_nY_l\xi_l \left\{-\frac{i\partial_{x_n}c(\xi')}{4(\xi_n-i)^2}
+h'(0)|\xi'|\left[\frac{11ic(\mathrm{d}x_n)}{8(\xi_n-i)}\right.\right.\nonumber\\
&\left.\left.-\frac{ic(\xi')+(1+4i)c(\mathrm{d}x_n)}{16(\xi_n-i)^2}-\frac{ic(\xi')-c(\mathrm{d}x_n)}{2(\xi_n-i)^3}\right]\right\} .\nonumber\\
\end{align}

We note that $i<n,~\int_{|\xi'|=1}\xi_{i_{1}}\xi_{i_{2}}\cdots\xi_{i_{2d+1}}\sigma(\xi')=0$,
so we omit some items that have no contribution for computing {\bf case a)~II)}.\\
Then there is the following formula
\begin{align}\label{c24}
&{\rm
tr} [\partial_{x_n}\pi^+_{\xi_n}\sigma_{-1}(\nabla_X^{S(TM)}\nabla_Y^{S(TM)}D^{-3})\times
\partial_{\xi_n}^2\sigma_{-3}(D^{-3})](x_0)\nonumber\\
&=-\sqrt{-1}\Sigma_{j,l=1}^{n-1}X_jY_l\xi_j\xi_lh'(0)\left[\frac{24i(5\xi_n^2-1)}{(\xi_n-i)^5(\xi_n+i)^4}
 +\frac{[2i(\xi_n-i)+26](5\xi_n^2-1)+24i(\xi_n^3-\xi_n)}{(\xi_n-i)^6(\xi_n+i)^4}\right.\nonumber\\
 &-\left.\frac{32(5\xi_n^2-1)+96(\xi_n^3-\xi_n)}{(\xi_n-i)^7(\xi_n+i)^4}\right]\nonumber\\
&+\sqrt{-1}X_nY_nh'(0)\left[\frac{46i(5\xi^2-1)-216(\xi_n^3-\xi_n)}{(\xi_n-i)^5(\xi_n+i)^4}
 -\frac{26(5\xi_n^2-1)+72i(\xi_n^3-\xi_n)}{(\xi_n-i)^6(\xi_n+i)^4}\right.\nonumber\\
 &-\left.\frac{16(5\xi_n^2-1)+48i(\xi_n^3-\xi_n)}{i(\xi_n-i)^7(\xi_n+i)^4}\right].\nonumber\\
\end{align}

Therefore, there is the following formula
\begin{align}\label{c25}
\Phi_2&=\frac{1}{2}\int_{|\xi'|=1}\int^{+\infty}_{-\infty}\bigg\{
-i\Sigma_{j,l=1}^{n-1}X_jY_l\xi_j\xi_lh'(0)\left[\frac{24i(5\xi_n^2-1)}{(\xi_n-i)^5(\xi_n+i)^4}\right.\nonumber\\
 &+\left.\frac{[2i(\xi_n-i)+26](5\xi_n^2-1)+24i(\xi_n^3-\xi_n)}{(\xi_n-i)^6(\xi_n+i)^4}
 -\frac{32(5\xi_n^2-1)+96(\xi_n^3-\xi_n)}{(\xi_n-i)^7(\xi_n+i)^4}\right]\nonumber\\
&+\sqrt{-1}X_nY_nh'(0)\left[\frac{46i(5\xi^2-1)-216(\xi_n^3-\xi_n)}{(\xi_n-i)^5(\xi_n+i)^4}
 -\frac{26(5\xi_n^2-1)+72i(\xi_n^3-\xi_n)}{(\xi_n-i)^6(\xi_n+i)^4}\right.\nonumber\\
 &-\left.\frac{16(5\xi_n^2-1)+48i(\xi_n^3-\xi_n)}{i(\xi_n-i)^7(\xi_n+i)^4}\right]
 \bigg\}d\xi_n\sigma(\xi')dx'\nonumber\\
 &=-\frac{1}{2}\sqrt{-1}\Sigma_{j,l=1}^{n-1}X_jY_lh'(0)\Omega_4\int_{\Gamma^{+}}
 \left[\frac{24i(5\xi_n^2-1)}{(\xi_n-i)^5(\xi_n+i)^4}
 +\frac{[2i(\xi_n-i)+26](5\xi_n^2-1)+24i(\xi_n^3-\xi_n)}{(\xi_n-i)^6(\xi_n+i)^4}\right.\nonumber\\
 &-\left.\frac{32(5\xi_n^2-1)+96(\xi_n^3-\xi_n)}{(\xi_n-i)^7(\xi_n+i)^4}\right]\xi_j\xi_ld\xi_{n}dx'\nonumber\\
 &+\frac{1}{2}\sqrt{-1}X_nY_nh'(0)\Omega_4\int_{\Gamma^{+}}
 \left[\frac{46i(5\xi^2-1)-216(\xi_n^3-\xi_n)}{(\xi_n-i)^5(\xi_n+i)^4}
 -\frac{26(5\xi_n^2-1)+72i(\xi_n^3-\xi_n)}{(\xi_n-i)^6(\xi_n+i)^4}\right.\nonumber\\
 &-\left.\frac{16(5\xi_n^2-1)+48i(\xi_n^3-\xi_n)}{i(\xi_n-i)^7(\xi_n+i)^4}\right]d\xi_{n}dx'\nonumber\\
&=-\frac{1}{2}\sqrt{-1}\Sigma_{j=1}^{n-1}X_jY_j\frac{4\pi}{3}h'(0)\Omega_4
\left\{ \frac{2\pi i}{4!}
\left[\frac{24i(5\xi_n^2-1)}{(\xi_n+i)^4}\right]^{(4)}\bigg|_{\xi_n=i}dx'\right.\nonumber\\
&+\frac{2\pi i}{5!}
\left[\frac{[2i(\xi_n-i)+26](5\xi_n^2-1)+24i(\xi_n^3-\xi_n)}{(\xi_n+i)^4}\right]^{(5)}\bigg|_{\xi_n=i}dx'\nonumber\\
 &\left.+\frac{2\pi i}{6!}
\left[\frac{32(5\xi_n^2-1)+96(\xi_n^3-\xi_n)}{(\xi_n+i)^4}\right]^{(6)}\bigg|_{\xi_n=i}dx'
\right\}\nonumber\\
 &+\frac{1}{2}\sqrt{-1}X_nY_nh'(0)\Omega_4\left\{
 \frac{2\pi i}{4!}
 \left[\frac{46i(5\xi^2-1)-216(\xi_n^3-\xi_n)}{(\xi_n+i)^4}\right]^{(4)}\bigg|_{\xi_n=i}dx'\right.\nonumber\\
 &\left.+\frac{2\pi i}{5!}
 \left[\frac{26(5\xi_n^2-1)+72i(\xi_n^3-\xi_n)}{(\xi_n+i)^4}\right]^{(5)}\bigg|_{\xi_n=i}dx'
 -\frac{2\pi i}{6!}
 \left[\frac{16(5\xi_n^2-1)+48i(\xi_n^3-\xi_n)}{(\xi_n+i)^4}\right]^{(6)}\bigg|_{\xi_n=i}dx'
\right\}\nonumber\\
&=\left[\left(-\frac{5925}{16}-\frac{877}{6}i\right)\pi\Sigma_{j=1}^{n-1}X_jY_j
+\left(311i\right)X_nY_n\right]h'(0)\pi\Omega_4dx',
\end{align}
where ${\rm \Omega_{4}}$ is the canonical volume of $S^{4}.$\\

\noindent  {\bf case a)~III)}~$r=-1,~l=-3,~j=|\alpha|=0,~k=1$.\\
\noindent By (\ref{c12}), we get
\begin{align}\label{c3}
\Phi_3&=-\frac{1}{2}\int_{|\xi'|=1}\int^{+\infty}_{-\infty}
{\rm trace} [\partial_{\xi_n}\pi^+_{\xi_n}\sigma_{-1}(\nabla_X^{S(TM)}\nabla_Y^{S(TM)}D^{-3})\times
\partial_{\xi_n}\partial_{x_n}\sigma_{-3}(D^{-3})](x_0)d\xi_n\sigma(\xi')dx'\nonumber\\
&=\frac{1}{2}\int_{|\xi'|=1}\int^{+\infty}_{-\infty}
{\rm trace} [\partial_{\xi_n}^2\pi^+_{\xi_n}\sigma_{-1}(\nabla_X^{S(TM)}\nabla_Y^{S(TM)}D^{-3})\times
\partial_{x_n}\sigma_{-3}(D^                                                                                   {-3})](x_0)d\xi_n\sigma(\xi')dx'.
\end{align}

\noindent By \cite{Ka}, we have
\begin{align}\label{c32}
\pi^+_{\xi_n}\sigma_{-1}(\nabla_X^{S(TM)}\nabla_Y^{S(TM)}D^{-3})
&=\left[\frac{c(\xi')}{4(\xi_n-i)}-\frac{ic(\xi')-c(\mathrm{d}x_n)}{4(\xi_n-i)^2}\right]
 \Sigma_{j,l=1}^{n-1}X_jY_l\xi_j\xi_l\nonumber\\
&+\left[\frac{c(\xi')-2ic(\mathrm{d}x_n)}{4(\xi_n-i)}+\frac{ic(\xi')-c(\mathrm{d}x_n)}{4(\xi_n-i)^2}\right]X_nY_n\nonumber\\
&+\left[\frac{c(\mathrm{d}x_n)}{4(\xi_n-i)}+\frac{c(\xi')+ic(\mathrm{d}x_n)}{4(\xi_n-i)^2}\right]\Sigma_{j=1}^{n-1}X_jY_n\xi_j\nonumber\\
&+\left[\frac{c(\mathrm{d}x_n)}{4(\xi_n-i)}+\frac{c(\xi')+ic(\mathrm{d}x_n)}{4(\xi_n-i)^2}\right]\Sigma_{l=1}^{n-1}X_nY_l\xi_l.
\end{align}

Then,
\begin{align}\label{c33}
\partial_{\xi_n}\pi^+_{\xi_n}\sigma_{-1}(\nabla_X^{S(TM)}\nabla_Y^{S(TM)}D^{-3})
&=\left[-\frac{ic(\xi')}{4(\xi_n-i)^2}+\frac{ic(\xi')-c(\mathrm{d}x_n)}{2(\xi_n-i)^3}\right]\Sigma_{j,l=1}^{n-1}X_jY_l\xi_j\xi_l\nonumber\\
&-\left[\frac{c(\xi')+2ic(\mathrm{d}x_n)}{4(\xi_n-i)^2}+\frac{ic(\xi')-c(\mathrm{d}x_n)}{2(\xi_n-i)^3}\right]X_nY_n\nonumber\\.
&-\left[\frac{c(\mathrm{d}x_n)}{4(\xi_n-i)^2}+\frac{c(\xi')+ic(\mathrm{d}x_n)}{2(\xi_n-i)^3}\right]\Sigma_{j=1}^{n-1}X_jY_n\xi_j\nonumber\\
&-\left[\frac{c(\mathrm{d}x_n)}{4(\xi_n-i)^2}+\frac{c(\xi')+ic(\mathrm{d}x_n)}{2(\xi_n-i)^3}\right]\Sigma_{l=1}^{n-1}X_nY_l\xi_l.
\end{align}
By (\ref{37}) and (\ref{c33}), we have
\begin{align}\label{c34}
&{\rm tr} [\partial_{\xi_n}\pi^+_{\xi_n}\sigma_{-1}(\nabla_X^{S(TM)}\nabla_Y^{S(TM)}D^{-3})\times
\partial_{\xi_n}\partial_{x_n}\sigma_{-3}(D^{-3})](x_0)\nonumber\\
&=\left[\frac{24i\xi_n}{(\xi_n-i)^2(1+\xi_n^2)^4}-\frac{4i\xi_n}{(\xi_n-i)^2(1+\xi_n^2)^3}
+\frac{48\xi_n-8i(1-5\xi_n^2)}{(\xi_n-i)^3(1+\xi_n^2)^4}-\frac{8\xi_n}{(\xi_n-i)^3(1+\xi_n^2)^3}\right]
h'(0)\sum_{j,l=1}^{n-1}X_jY_l\xi_j\xi_l\nonumber\\
&+\left[\frac{24i\xi_n+8(1_5\xi_n^2)}{(\xi_n-i)^2(1+\xi_n^2)^4}-\frac{4i\xi_n}{(\xi_n-i)^2(1+\xi_n^2)^3}
-\frac{48\xi_n-8i(1-5\xi_n^2)}{(\xi_n-i)^3(1+\xi_n^2)^4}+\frac{8\xi_n}{(\xi_n-i)^3(1+\xi_n^2)^3}\right]
h'(0)X_nY_n.\nonumber\\
\end{align}

Therefore, we get the following formula
\begin{align}\label{c35}
\Phi_3&=\frac{1}{2}\int_{|\xi'|=1}\int^{+\infty}_{-\infty}
\bigg(\left[\frac{24i\xi_n}{(\xi_n-i)^2(1+\xi_n^2)^4}-\frac{4i\xi_n}{(\xi_n-i)^2(1+\xi_n^2)^3}
   +\frac{48\xi_n-8i(1-5\xi_n^2)}{(\xi_n-i)^3(1+\xi_n^2)^4}\right.\nonumber\\
&-\left.\frac{8\xi_n}{(\xi_n-i)^3(1+\xi_n^2)^3}\right]
   h'(0)\sum_{j,l=1}^{n-1}X_jY_l\xi_j\xi_l\nonumber\\
&+\left[\frac{24i\xi_n+8(1_5\xi_n^2)}{(\xi_n-i)^2(1+\xi_n^2)^4}-\frac{4i\xi_n}{(\xi_n-i)^2(1+\xi_n^2)^3}
  -\frac{48\xi_n-8i(1-5\xi_n^2)}{(\xi_n-i)^3(1+\xi_n^2)^4}\right.\nonumber\\
&+\left.\frac{8\xi_n}{(\xi_n-i)^3(1+\xi_n^2)^3}\right]
   h'(0)X_nY_n bigg)d\xi_n\sigma(\xi')dx'\nonumber\\
&=-\frac{1}{2}\Sigma_{j,l=1}^{n-1}X_jY_lh'(0)\Omega_4\left\{
\frac{2\pi i}{5!}\left[\frac{24\xi_n}{(\xi_n+i)^4}\right]^{(5)}\bigg|_{\xi_n=i}
-\frac{2\pi i}{4!}\left[\frac{4\xi_n}{(\xi_n+i)^3}\right]^{(4)}\bigg|_{\xi_n=i}\right.\nonumber\\
&+\left.\frac{2\pi i}{6!}
  \left[\frac{48\xi_n-8i(1-5\xi_n^2)}{(\xi_n+i)^4}\right]^{(6)}\bigg|_{\xi_n=i}
-\frac{2\pi i}{5!}\left[\frac{8\xi_n}{(\xi_n+i)^3}\right]^{(5)}\bigg|_{\xi_n=i}
\right\}dx'\nonumber\\
&-\frac{1}{2}X_nY_nh'(0)\Omega_4\left\{
\frac{2\pi i}{5!}
 \left[\frac{24\xi_n+8(1-5\xi_n^2)}{(\xi_n+i)^4}\right]^{(5)}\bigg|_{\xi_n=i}
-\frac{2\pi i}{4!}
 \left[\frac{4i\xi_n}{(\xi_n+i)^3}\right]^{(4)}\bigg|_{\xi_n=i}\right.\nonumber\\
&-\left.\frac{2\pi i}{6!}
 \left[\frac{48\xi_n-8i(1-5\xi_n^2)}{(\xi_n+i)^4}\right]^{(6)}\bigg|_{\xi_n=i}
+\frac{2\pi i}{5!}\left[\frac{8xi_n}{(\xi_n+i)^3}\right]^{(5)}\bigg|_{\xi_n=i}
\right\}dx'\nonumber\\
&=\left[\left(\frac{49}{3}+86i\right)\pi\Sigma_{j=1}^{n-1}X_jY_j-\left(\frac{17}{4}-68i\right)X_nY_n\right]h'(0)\pi\Omega_3dx'.
\end{align}

\noindent  {\bf case b)}~$r=-1,~l=-4,~k=j=|\alpha|=0$.\\
\noindent By (\ref{c12}), we get
\begin{align}\label{c4}
\Phi_4&=-i\int_{|\xi'|=1}\int^{+\infty}_{-\infty}{\rm trace} [\pi^+_{\xi_n}\sigma_{-1}(\nabla_X^{S(TM)}\nabla_Y^{S(TM)}D^{-3})\times
\partial_{\xi_n}\sigma_{-4}(D^{-3})](x_0)d\xi_n\sigma(\xi')dx'\nonumber\\
&=i\int_{|\xi'|=1}\int^{+\infty}_{-\infty}{\rm trace} [\partial_{\xi_n}\pi^+_{\xi_n}\sigma_{-1}(\nabla_X^{S(TM)}\nabla_Y^{S(TM)}D^{-3})\times
\sigma_{-4}(D^{-3})](x_0)d\xi_n\sigma(\xi')dx'.
\end{align}

We note that $i<n,~\int_{|\xi'|=1}\xi_{i_{1}}\xi_{i_{2}}\cdots\xi_{i_{2d+1}}\sigma(\xi')=0$,
so we omit some items that have no contribution for computing {\bf case b)}.\\
By (\ref{62}) and (\ref{c33}), we have

\begin{align}\label{c42}
&{\rm tr}[\partial_{\xi_n}\pi^+_{\xi_n}\sigma_{-1}(\nabla_X^{S(TM)}\nabla_Y^{S(TM)}D^{-3})\times
\sigma_{-3}(D^{-3})](x_0)\nonumber\\
&=h'(0)\left[\frac{(-11+3i)\xi_n(1+\xi_n^2)-16i\xi_n}{(\xi_n-i)^2(1+\xi_n^2)^4}\right.\nonumber\\
&-\left.\sqrt{-1}\frac{(-22+6i)\xi_n(1+\xi_n^2)-32i\xi_n+(26-7i)(1+\xi_n^2)-32-15i\xi_n^2(1+\xi_n^2)}{(\xi_n-i)^3(\xi_n+i)^4}\right]
\Sigma_{j,l=1}^{n-1}X_jY_l\xi_j\xi_l\nonumber\\
&-\sqrt{-1}\left[\frac{(-11-3i)\xi_n(1+\xi_n^2)+16\xi_n+(48i-14)(1+\xi_n^2)-32i+30\xi^2(1+\xi_n^2)}{(\xi_n-i)^6(\xi_n+i)^4}\right.\nonumber\\
&+\left.\frac{(-22+6i)\xi_n(1+\xi_n^2)-32i\xi_n+(26-7i)(1+\xi_n^2)-32-15i\xi_n^2(1+\xi_n^2)}{(\xi_n-i)^3(\xi_n+i)^4}\right]
h'(0)X_nY_n.\nonumber\\
\end{align}

Therefore, we get
\begin{align}\label{c43}
\Phi_4&=i\int_{|\xi'|=1}\int^{+\infty}_{-\infty}
 h'(0)\left[\frac{(-11+3i)\xi_n(1+\xi_n^2)-16i\xi_n}{(\xi_n-i)^2(1+\xi_n^2)^4}\right.\nonumber\\
 &-\left.\sqrt{-1}\frac{(-22+6i)\xi_n(1+\xi_n^2)-32i\xi_n+(26-7i)(1+\xi_n^2)-32-15i\xi_n^2(1+\xi_n^2)}{(\xi_n-i)^3(1+\xi_n^2)^4}\right]
 \Sigma_{j,l=1}^{n-1}X_jY_l\xi_j\xi_l\nonumber\\
 &-\sqrt{-1}\left[\frac{(-11-3i)\xi_n(1+\xi_n^2)+16\xi_n+(48i-14)(1+\xi_n^2)-32i+30\xi^2(1+\xi_n^2)}{(\xi_n-i)^2(1+\xi_n^2)^4}\right.\nonumber\\
 &+\left.\frac{(-22+6i)\xi_n(1+\xi_n^2)-32i\xi_n+(26-7i)(1+\xi_n^2)-32-15i\xi_n^2(1+\xi_n^2)}{(\xi_n-i)^3(1+\xi_n^2)^4}\right]
 h'(0)X_nY_nd\xi_n\sigma(\xi')dx'\nonumber\\
&=\Sigma_{j,l=1}^{n-1}X_jY_lh'(0)\Omega_4\int_{\Gamma^{+}}\left[
 \frac{(-11+3i)\xi_n}{(\xi_n-i)^2(1+\xi_n^2)^3}
 -\frac{16i\xi_n}{(\xi_n-i)^2(1+\xi_n^2)^4}
 -\frac{(-22+6i)i\xi_n}{(\xi_n-i)^3(\xi_n+i)^3}\right.\nonumber\\
 &-\left.\frac{32\xi_n}{(\xi_n-i)^3(1+\xi_n^2)^4}
 -\frac{(26-7i)i}{(\xi_n-i)^3(1+\xi_n^2)^3}
 +\frac{32i}{(\xi_n-i)^3(1+\xi_n^2)^4}
 -\frac{15\xi_n^2}{(\xi_n-i)^3(1+\xi_n^2)^3}
 \right]\xi_j\xi_ld\xi_{n}dx'\nonumber\\
 &+X_nY_nh'(0)\Omega_4\int_{\Gamma^{+}}\left[
 \frac{16\xi_n}{(\xi_n-i)^6(\xi_n+i)^4}
 -\frac{(3+11i)\xi_n}{(\xi_n-i)^5(\xi_n+i)^3}
 +\frac{48i-14}{(\xi_n-i)^6(\xi_n+i)^4}\right.\nonumber\\
 &-\frac{32i}{(\xi_n-i)^5(\xi_n+i)^3}
 +\frac{30\xi_n^2}{(\xi_n-i)^5(\xi_n+i)^3}
 +\frac{(-22+6i)i\xi_n}{(\xi_n-i)^3(\xi_n+i)^3}
 +\frac{32\xi_n}{(\xi_n-i)^3(1+\xi_n^2)^4}\nonumber\\
 &+\left.\frac{(26-7i)i}{(\xi_n-i)^3(1+\xi_n^2)^3}
 -\frac{32i}{(\xi_n-i)^3(1+\xi_n^2)^4}
 +\frac{15\xi_n^2}{(\xi_n-i)^3(1+\xi_n^2)^3}
 \right]d\xi_{n}dx'\nonumber\\
&=\Sigma_{j,l=1}^{n-1}X_jY_lh'(0)\pi^2\Omega_4\left\{
 \frac{2\pi i}{4!}\left[\frac{(-11+3i)\xi_n}{(\xi_n+i)^3}\right]^{(4)}\bigg|_{\xi_n=i}dx'
 -\frac{2\pi i}{5!}\left[\frac{16i\xi_n}{(\xi_n+i)^4}\right]^{(5)}\bigg|_{\xi_n=i}dx'\right.\nonumber\\
 &-\frac{2\pi i}{5!}\left[\frac{(-22+6i)i\xi_n}{(\xi_n+i)^3}\right]^{(5)}\bigg|_{\xi_n=i}dx'
 -\frac{2\pi i}{6!}\left[\frac{32\xi_n}{(\xi_n+i)^4}\right]^{(6)}\bigg|_{\xi_n=i}dx'
 -\frac{2\pi i}{5!}\left[\frac{(26-7i)i}{(\xi_n+i)^3}\right]^{(5)}\bigg|_{\xi_n=i}dx'\nonumber\\
 &\left.+\frac{2\pi i}{6!}\left[\frac{32i}{(\xi_n+i)^4}\right]^{(6)}\bigg|_{\xi_n=i}dx'
 -\frac{2\pi i}{5!}\left[\frac{15\xi_n^2}{(\xi_n+i)^3}\right]^{(5)}\bigg|_{\xi_n=i}dx'
 \right\}\nonumber\\
 &+X_nY_nh'(0)\Omega_4\left\{\
 \frac{2\pi i}{5!}\left[\frac{16\xi_n}{(\xi_n+i)^4}\right]^{(5)}\bigg|_{\xi_n=i}dx'
 -\frac{2\pi i}{4!}\left[\frac{(3+11i)\xi_n}{(\xi_n+i)^3}\right]^{(4)}\bigg|_{\xi_n=i}dx'\right.\nonumber\\
 &-\frac{2\pi i}{5!}\left[\frac{48i-14}{(\xi_n+i)^4}\right]^{(5)}\bigg|_{\xi_n=i}dx'
 -\frac{32i\times2\pi i}{4!}\left[\frac{1}{(\xi_n+i)^3}\right]^{(4)}\bigg|_{\xi_n=i}dx'
 +\frac{2\pi i}{4!}\left[\frac{30\xi_n^2}{(\xi_n+i)^3}\right]^{(4)}\bigg|_{\xi_n=i}dx'\nonumber\\
 &+\frac{2\pi i}{5!}\left[\frac{(-22+6i)i\xi_n}{(\xi_n+i)^3}\right]^{(5)}\bigg|_{\xi_n=i}dx'
 +\frac{2\pi i}{6!}\left[\frac{32\xi_n}{(\xi_n+i)^4}\right]^{(6)}\bigg|_{\xi_n=i}dx'
 +\frac{2\pi i}{5!}\left[\frac{(26-7i)i}{(\xi_n+i)^3}\right]^{(5)}\bigg|_{\xi_n=i}dx'\nonumber\\
 &\left.-\frac{2\pi i}{6!}\left[\frac{32i}{(\xi_n+i)^4}\right]^{(6)}\bigg|_{\xi_n=i}dx'
 +\frac{2\pi i}{5!}\left[\frac{15\xi_n^2}{(\xi_n+i)^3}\right]^{(5)}\bigg|_{\xi_n=i}dx'
dx'\right\}\nonumber\\
&=\left[-\left(\frac{53}{3}+\frac{41}{6}i\right)\pi\Sigma_{j=1}^{n-1}X_jY_j
 +\left(\frac{173i}{4}+36\right)X_nY_n\right]h'(0)\pi\Omega_3dx'.
\end{align}

\noindent {\bf  case c)}~$r=-2,~\ell=-3,~k=j=|\alpha|=0$.\\
By (\ref{c12}), we get
\begin{align}\label{c5}
\Phi_5&=-\int_{|\xi'|=1}\int^{+\infty}_{-\infty}{\rm trace} [\pi^+_{\xi_n}\sigma_{-2}(\nabla_X^{S(TM)}\nabla_Y^{S(TM)}D^{-3})\times
\partial_{\xi_n}\sigma_{-3}(D^{-3})](x_0)d\xi_n\sigma(\xi')dx'.
\end{align}

By \cite{Ka}, we have
\begin{align}\label{c53}
\sigma_{-2}(\nabla_X^{S(TM)}\nabla_Y^{S(TM)}D^{-3})
&=-\sum_j A(Y)X_j\xi_j\frac{c(\xi)}{|\xi|^{4}}-\sum_l A(X)Y_l\xi_l\frac{c(\xi)}{|\xi|^{4}}
-\sum_{j,l=1}^{n-1}X_j\frac{\partial_{Y_l}}{\partial_{X_j}}\xi_l\frac{c(\xi)}{|\xi|^{4}}\nonumber\\
&+\sum_{j,l=1}^nX_jY_l\xi_j\xi_l\frac{1}{(\xi_n^2+1)^4}
\left\{\left(\frac{11}{2}\xi_n(1+\xi_n^2)+8i\xi_n\right)h'(0)c(\xi')\right.\nonumber\\
&+\left[-2i+6i\xi_n^2-\frac{7}{4}(1+\xi_n^2)
 +\frac{15}{4}\xi_n^2(1+\xi^2_n)\right]h'(0)c(\mathrm{d}x_n) \nonumber\\
&\left.-3i\xi_n(1+\xi^2_n)\partial_{x_n}c(\xi')
 +i(1+\xi^2_n)c(\xi')c(\mathrm{d}x_n)\partial_{x_n}c(\xi')\right\}.
\end{align}

We note that $i<n,~\int_{|\xi'|=1}\xi_{i_{1}}\xi_{i_{2}}\cdots\xi_{i_{2d+1}}\sigma(\xi')=0$,
so we omit some items that have no contribution for computing {\bf case c)}.\\

By (\ref{c53}) and (\ref{43}), we have
\begin{align}\label{c54}
&{\rm tr}[\pi^+_{\xi_n}\sigma_{-2}(\nabla_X^{S(TM)}\nabla_Y^{S(TM)}D^{-3})\times
\partial_{\xi_n}\sigma_{-3}(D^{-3})](x_0)|_{|\xi'|=1}\nonumber\\
&=\left[\frac{1-9i}{(\xi_n-i)^4(\xi_n+i)^3}-\frac{6-27i}{4(\xi_n-i)^3(\xi_n+i)^2}
  -\frac{27+58i}{(\xi_n-i)^5(\xi_n+i)^3}+\frac{(1+5i)\xi_n}{(\xi_n-i)^5(\xi_n+i)^3}\right.\nonumber\\
&+\frac{81+174i}{4(\xi_n-i)^4(\xi_n+i)^2}-\frac{10-2i}{(\xi_n-i)^6(\xi_n+i)^3}
 +\frac{(14-i)\xi_n}{(\xi_n-i)^6(\xi_n+i)^3}+\frac{16-3i}{2(\xi_n-i)^5(\xi_n+i)^2}\nonumber\\
&+\left.\frac{4i-12}{(\xi_n-i)^7(\xi_n+i)^3}-\frac{16i\xi_n}{(\xi_n-i)^7(\xi_n+i)^3}
+\frac{9-3i}{(\xi_n-i)^6(\xi_n+i)^2}\right]
 h'(0)\sum_{j,l=1}^{n-1}X_jY_l\xi_j\xi_l\nonumber\\
&+\left[-\frac{72+25i}{4(\xi_n-i)^4(\xi_n+i)^3}+\frac{54+8i}{(\xi_n-i)^3(\xi_n+i)^2}
-\frac{7+43i}{(\xi_n-i)^5(\xi_n+i)^3}+\frac{(44-9i)\xi_n}{2(\xi_n-i)^5(\xi_n+i)^3}\right.\nonumber\\
&+\frac{21+129i}{4(\xi_n-i)^4(\xi_n+i)^2}+\frac{22-26i}{(\xi_n-i)^6(\xi_n+i)^3}
+\frac{3\xi_n}{(\xi_n-i)^6(\xi_n+i)^3}-\frac{33-39i}{2(\xi_n-i)^6(\xi_n+i)^2}\nonumber\\
&+\left.\frac{2i}{(\xi_n-i)^7(\xi_n+i)^3}+\frac{16\xi_n}{(\xi_n-i)^7(\xi_n+i)^3}
 -\frac{3i}{2(\xi_n-i)^6(\xi_n+i)^2}\right]h'(0)X_nY_n.
\end{align}

So we have
\begin{align}\label{c55}
\Phi_5&=-i\int_{|\xi'|=1}\int^{+\infty}_{-\infty}{\rm tr}[\pi^+_{\xi_n}\sigma_{1}(\nabla_X^{S(TM)}\nabla_Y^{S(TM)}D^{-1})\times
\partial_{\xi_n}\sigma_{-4}(D^{-3})](x_0)d\xi_n\sigma(\xi')dx'\nonumber\\
&=\Sigma_{j,l=1}^{n-1}X_jY_lh'(0)\Omega_4\int_{\Gamma^{+}}
 \left[\frac{1-9i}{(\xi_n-i)^4(\xi_n+i)^3}-\frac{6-27i}{4(\xi_n-i)^3(\xi_n+i)^2}
  -\frac{27+58i}{(\xi_n-i)^5(\xi_n+i)^3}\right.\nonumber\\
  &+\frac{(1+5i)\xi_n}{(\xi_n-i)^5(\xi_n+i)^3}
   +\frac{81+174i}{4(\xi_n-i)^4(\xi_n+i)^2}-\frac{10-2i}{(\xi_n-i)^6(\xi_n+i)^3}
   +\frac{(14-i)\xi_n}{(\xi_n-i)^6(\xi_n+i)^3}\nonumber\\
  &+\left.\frac{16-3i}{2(\xi_n-i)^5(\xi_n+i)^2}
   +\frac{4i-12}{(\xi_n-i)^7(\xi_n+i)^3}-\frac{16i\xi_n}{(\xi_n-i)^7(\xi_n+i)^3}
   +\frac{9-3i}{(\xi_n-i)^6(\xi_n+i)^2}\right]\xi_j\xi_ld\xi_{n}dx'\nonumber\\
&+X_nY_nh'(0)\Omega_4\int_{\Gamma^{+}}
\left[-\frac{72+25i}{4(\xi_n-i)^4(\xi_n+i)^3}+\frac{54+8i}{(\xi_n-i)^3(\xi_n+i)^2}
 -\frac{7+43i}{(\xi_n-i)^5(\xi_n+i)^3}\right.\nonumber\\
 &+\frac{(44-9i)\xi_n}{2(\xi_n-i)^5(\xi_n+i)^3}
 +\frac{21+129i}{4(\xi_n-i)^4(\xi_n+i)^2}+\frac{22-26i}{(\xi_n-i)^6(\xi_n+i)^3}
 +\frac{3\xi_n}{(\xi_n-i)^6(\xi_n+i)^3}\nonumber\\
 &-\left.\frac{33-39i}{2(\xi_n-i)^6(\xi_n+i)^2}
 +\frac{2i}{(\xi_n-i)^7(\xi_n+i)^3}+\frac{16\xi_n}{(\xi_n-i)^7(\xi_n+i)^3}
 -\frac{3i}{2(\xi_n-i)^6(\xi_n+i)^2}\right]d\xi_{n}dx'\nonumber\\
&=\Sigma_{j,l=1}^{n-1}X_jY_lh'(0)\pi\Omega_4\left\{
 \frac{2\pi i}{3!}\left[\frac{1-9i}{(\xi_n+i)^3}\right]^{(3)}\bigg|_{\xi_n=i}dx'
 -\frac{2\pi i}{2!}\left[\frac{6-27i}{4(\xi_n+i)^2}\right]^{(2)}\bigg|_{\xi_n=i}dx'\right.\nonumber\\
 &-\frac{2\pi i}{4!}\left[\frac{27+58i}{(\xi_n+i)^3}\right]^{(4)}\bigg|_{\xi_n=i}dx'
 +\frac{2\pi i}{4!}\left[\frac{(1+5i)\xi_n}{(\xi_n+i)^3}\right]^{(4)}\bigg|_{\xi_n=i}dx'
 +\frac{2\pi i}{3!}\left[\frac{81+174i}{4(\xi_n+i)^2}\right]^{(3)}\bigg|_{\xi_n=i}dx'\nonumber\\
 &-\frac{2\pi i}{5!}\left[\frac{10-2i}{(\xi_n+i)^3}\right]^{(5)}\bigg|_{\xi_n=i}dx'
  +\frac{2\pi i}{5!}\left[\frac{(14-i)\xi_n}{(\xi_n+i)^3}\right]^{(5)}\bigg|_{\xi_n=i}dx'
  +\frac{2\pi i}{4!}\left[\frac{16-3i}{2(\xi_n+i)^2}\right]^{(4)}\bigg|_{\xi_n=i}dx'\nonumber\\
 &+\left.\frac{2\pi i}{6!}\left[\frac{4i-12}{(\xi_n+i)^3}\right]^{(6)}\bigg|_{\xi_n=i}dx'
 -\frac{2\pi i}{6!}\left[\frac{16i\xi_n}{(\xi_n+i)^3}\right]^{(6)}\bigg|_{\xi_n=i}dx'
  +\frac{2\pi i}{5!}\left[\frac{9-3i}{(\xi_n+i)^2}\right]^{(5)}\bigg|_{\xi_n=i}dx'
\right\}\nonumber\\
&+X_nY_nh'(0)\Omega_4\left\{
\frac{2\pi i}{3!}\left[-\frac{72+25i}{4(\xi_n+i)^3}\right]^{(3)}\bigg|_{\xi_n=i}dx'
+\frac{2\pi i}{2!}\left[\frac{54+8i}{(\xi_n+i)^2}\right]^{(2)}\bigg|_{\xi_n=i}dx'\right.\nonumber\\
 &-\frac{2\pi i}{4!}\left[\frac{7+43i}{(\xi_n+i)^3}\right]^{(4)}\bigg|_{\xi_n=i}dx'
 +\frac{2\pi i}{4!}\left[\frac{(44-9i)\xi_n}{2(\xi_n+i)^3}\right]^{(4)}\bigg|_{\xi_n=i}dx'
 +\frac{2\pi i}{3!}\left[\frac{21+129i}{4(\xi_n+i)^2}\right]^{(3)}\bigg|_{\xi_n=i}dx'\nonumber\\
 &+\frac{2\pi i}{5!}\left[\frac{22-26i}{(\xi_n+i)^3}\right]^{(5)}\bigg|_{\xi_n=i}dx'
 +\frac{2\pi i}{5!}\left[\frac{3\xi_n}{(\xi_n+i)^3}\right]^{(5)}\bigg|_{\xi_n=i}dx'
 -\frac{2\pi i}{5!}\left[\frac{33-39i}{2(\xi_n+i)^2}\right]^{(5)}\bigg|_{\xi_n=i}dx'\nonumber\\
 &+\left.\frac{2\pi i}{6!}\left[\frac{2i}{(\xi_n+i)^3}\right]^{(6)}\bigg|_{\xi_n=i}dx'
 +\frac{2\pi i}{6!}\left[\frac{16\xi_n}{(\xi_n+i)^3}\right]^{(6)}\bigg|_{\xi_n=i}dx'
 -\frac{2\pi i}{5!}\left[\frac{3i}{2(\xi_n+i)^2}\right]^{(5)}\bigg|_{\xi_n=i}dx'
 \right\}\nonumber\\
&=\left[\left(-\frac{775}{6}+113i\right)\sum_{j=1}^{n-1}X_jY_j\pi+\left(\frac{171}{2}+\frac{369}{8}i\right)X_nY_n\right]h'(0)\pi\Omega_4dx'.
\end{align}

Let $X=X^T+X_n\partial_n,~Y=Y^T+Y_n\partial_n,$ then we have $\sum_{j=1}^{n-1}X_jY_j=g(X^T,Y^T).$ Now $\Phi$ is the sum of the cases (a), (b) and (c). Therefore, we get
\begin{align}\label{c01}
\Phi=\sum_{i=1}^5\Phi_i=\left[\left(\frac{475}{4}+\frac{3747}{8}i\right)X_nY_n+\left(-\frac{8013}{16}+46i\right)\pi g(X^T,Y^T)\right]h'(0)\pi\Omega_4dx'.
\end{align}
Then, by (\ref{c1})-(\ref{c01}), we obtain following theorem

\begin{thm}\label{thmc1}
Let $M$ be a $6$-dimensional oriented
compact spin manifold with boundary $\partial M$ and the metric
$g^{M}$ be defined as Sect.(\ref{section:3}), then we get the following equality:
\begin{align}
\label{c02}
&\widetilde{{\rm Wres}}[\pi^+(\nabla_X^{S(TM)}\nabla_Y^{S(TM)}D^{-3})\circ\pi^+(D^{-3})]\nonumber\\
&=\frac{4\pi^3}{3}\int_{M}[Ric(X,Y)-\frac{1}{2}sg(X,Y)]d{\rm Vol_{M}}+\pi^3\int_{M}sg(X,Y)d{\rm Vol_{M}}\nonumber\\
&+\int_{\partial M}\left[\left(\frac{475}{4}+\frac{3747}{8}i\right)X_nY_n+\left(-\frac{8013}{16}+46i\right)\pi g(X^T,Y^T)\right]h'(0)\pi\Omega_4d{\rm Vol_{M}}.
\end{align}
\end{thm}

\section*{Declarations}
\begin{itemize}
\item Ethics approval and consent to participate\\
Not applicable.
\item Consent for publication\\
Not applicable.
\item Availability of data and materials\\
The authors confirm that the data supporting the findings of this study
are available within the article.
\item Competing interests\\
The authors declare no conflict of interest.
\item Funding\\
This research was funded by National Natural Science Foundation of China:
No.11771070.
\item Authors' contributions\\
All authors contributed to the study conception and design. Material
preparation, data collection and analysis were performed by YY and TW. The first draft of the
manuscript was written by YY and all authors commented on previous versions of the manuscript.
All authors read and approved the final manuscript.
\item Acknowledgements\\
This work was supported by  NSFC No.11771070. The authors thank the referee for his (or her) careful reading and helpful comments..
\end{itemize}

\end{document}